\documentclass[10pt]{article} 
\usepackage{fullpage}
\usepackage[T1]{fontenc} 
\usepackage[latin1]{inputenc}
\usepackage[english]{babel}
\usepackage{selinput}
\SelectInputMappings{adieresis={ä},germandbls={ß}}

 \usepackage{hyperref}

\usepackage{amsmath} 
\usepackage{amssymb} 
\usepackage{amsthm} 
\usepackage{enumerate}

\usepackage{graphicx} 
\usepackage{epstopdf}
\usepackage{subfig}
\usepackage[framed,numbered]{matlab-prettifier}
\usepackage{algorithm} 
\usepackage[noend]{algpseudocode}

\graphicspath{{./Figures/}}

\usepackage{multirow}
\usepackage{fancyvrb}

\newcommand{\R}{\mathbb{R}}

\newcommand{\set}[1]{\left\lbrace #1 \right\rbrace}

\usepackage{ mathrsfs }
\renewcommand{\L}{\mathcal{L}}

\usepackage[numbers,sort]{natbib}

\newcommand{\rev}[1]{\textcolor{black}{#1}}

\title{Fast global spectral methods for three-dimensional partial differential equations}
\usepackage{authblk}
\author[1]{Christoph Strössner}
\author[2]{Daniel Kressner}
\affil[1]{Institute   of   Mathematics, EPF  Lausanne,  Switzerland,  \url{christoph.stroessner@epfl.ch}.}
\affil[2]{Institute   of   Mathematics, EPF  Lausanne,  Switzerland,  \url{daniel.kressner@epfl.ch}.}
\date{March 17, 2022}

\begin{document}
\maketitle

\begin{abstract} 
Global spectral methods offer the potential to compute solutions of partial differential equations numerically to very high accuracy. 
In this work, we develop a novel global spectral method for linear partial differential equations on cubes by extending ideas of Chebop2 [Townsend and Olver, \emph{J. Comput. Phys.}, 299 (2015)] to the three-dimensional setting utilizing expansions in tensorized polynomial bases.
Solving the discretized PDE involves a linear system that can be recast as a linear tensor equation. Under suitable additional assumptions, the structure of these equations admits for an efficient solution via the blocked recursive solver [Chen and Kressner, \emph{Numer. Algorithms}, 84 (2020)].
In the general case, when these assumptions are not satisfied, this solver is used as a preconditioner to speed up computations.
\end{abstract}

\section{Introduction}
This work is concerned with the solution of linear partial differential equations (PDEs) on cubes using global spectral methods~\cite{Fornberg96,Gottlieb77,Townsend15,Trefethen01}. 
The solution is approximated globally on the whole domain in terms of a truncated (tensorized) Chebyshev series~\cite{Zhao94,Mason02,Trefethen13}.  
For one- and two-dimensional rectangular domains this approach led to the development of the solvers Chebop~\cite{Driscoll08,Olver13,Driscoll16} and Chebop2~\cite{Townsend15} contained in the Chebfun package~\cite{Driscoll14} for computing numerically with functions. 
Many models used in medicine~\cite{Hillen09,Liu11,Swanson03}, engineering~\cite{Bazilevs11,Ebling09,Popov00} and geosciences~\cite{Igel99,McBride06,Zhdanov06} involve PDEs on three-dimensional domains. 
These domains are rarely cubes, but they can often be mapped onto cubes~\rev{\cite{Fortunato21}}.
Whilst the efficient approximation of trivariate functions on cubes by combining tensorized Chebyshev interpolation and low-rank approximations has been studied in Chebfun3~\cite{Hashemi17} and Chebfun3F~\cite{Stroessner21}, there is no generalization of Chebop2 for three-dimensional PDEs on cubes. 
So far, the ideas of Chebop2 have been extended to triangular domains~\cite{Olver19} and to disks~\cite{Wilber17,Olver20}. 
For the special case of the Poisson\rev{'s} equation with homogeneous Dirichlet boundary conditions, extensions to three-dimensional spheres~\cite{Townsend16b}, cylinders and unit cubes~\cite{Fortunato19} have been developed.

Throughout this work, we study linear PDEs of the form
\begin{align}\label{eq:SystemIntroduction}
    \L u = f \quad \textrm{on } [-1,1]^3,
\end{align}
complemented with linear boundary conditions. 
We approximate $u$ and $f$ in terms of truncated expansions with \rev{tensorized} polynomial basis functions.
These expansions are represented by coefficient tensors of order 3. 
Global spectral methods approximate the operator $\L$ by a mapping \rev{mimicking} the impact of \rev{applying} $\L$ on the level of the coefficient tensors.
We compute such a mapping by expressing the three-dimensional differential operator as combination of one-dimensional linear differential operators, which we determine using a \rev{canonical polyadic (CP) decomposition}~\cite{Kolda09}. 
By representing $u$ in a Chebyshev basis and $f$ in an ultraspherical basis, we can express the action of these one-dimensional operators on the coefficients in terms of sparse and well-conditioned differentiation and multiplication matrices~\cite{Olver13}.
The combination of applying these matrices yields the desired mapping. 
Solving the PDE corresponds to inverting this mapping, which can be seen as solving a tensor-valued linear system.
Discretizing the boundary conditions provides additional linear constraints under which inverting the linear system has a unique solution. 
Using substitution we obtain an unconstrained tensor-valued linear system uniquely determining a subtensor, from \rev{which} we recover the full coefficient tensor of the solution.
\rev{We would like to point out that our discretization approach can be applied to a wide variety of PDEs, but it does not necessarily preserve certain desirable properties of the differential operator.}

It turns out that our discretization of the PDE~\eqref{eq:SystemIntroduction} has Kronecker structure, which can often be exploited. For instance, Poisson's equation with homogeneous Dirichlet boundary conditions leads to a Laplace-like equation~\cite{Ballani13b,Chen12,Kressner09,Muravleva19}. Laplace-like equations can be solved efficiently using the blocked recursive algorithm in~\cite{Chen19}. This algorithm is asymptotically slower than the  nested alternating direction implicit method in~\cite{Fortunato19}, but our numerical experiments in Section~\ref{sec:NumericalExperiments} demonstrate that the recursive blocked algorithm is significantly faster in practice. Even if the PDE of interest does not immediately lead to  a structure suitable for the fast solver, we can often apply the blocked recursive algorithm as \rev{a} preconditioner for GMRES. 
 
The algorithmic ideas presented in this work, can \rev{be extended from solving} stationary linear PDEs, such as the Helmholtz equation and (convection) diffusion problems, to \rev{solving} time-dependent PDEs and PDE eigenvalue problems of the form
\[
\frac{\partial}{\partial t} u = \L u \quad \text{and} \quad \L u = \lambda u,
\]
which we solve using implicit Euler and inverse iteration methods\rev{,} respectively.
\rev{We expect to obtain accurate approximations of the solution, when both the PDE coefficients and the solution are sufficiently smooth to be well approximated by truncated expansions with tensorized polynomial basis functions. For non-smooth solutions, our global spectral method might lead to inaccurate approximations.}

The global spectral method can be used to compute solutions numerically to very high accuracy. It is not to be confused with so\rev{-}called spectral element methods~\cite{Fortunato21,Hesthaven07,Karniadakis05} and $p$- and $hp$-finite element methods~\cite{Babuska90,Olver20,Xu18}\rev{. In those methods, $u$ is not approximated globally by a truncated series expansion. Instead, $u$ is written as sum of (locally supported) functions, each of which is approximated individually by a truncated series expansion.} There also exist solvers relying on using domain decompositions in combination with truncated series expansions~\cite{Hao16,Pfeiffer03}.
We want to emphasize that the methods presented in this work can be used as local solver, when the elements/subdomains can be mapped to cubes.

The remainder of this paper is structured as follows. 
In Section~\ref{sec:problemsetting}, we define linear differential operators and the approximation format. 
In Section~\ref{sec:Discretization}, we derive the mapping on the level of the coefficient tensor for linear differential operators. 
The discretization of PDEs and the efficient solution of the resulting tensor-valued linear system is discussed in Section~\ref{sec:SpectralMethod}.
In Section~\ref{sec:NumericalExperiments}, we apply our global spectral method to solve stationary PDEs, parabolic PDEs and PDE eigenvalue problems numerically to very high accuracy.

\section{Problem setting}\label{sec:problemsetting}

\subsection{Structure of a linear differential operator}
A linear partial differential operator $\L$ on the domain $[-1,1]\times[-1,1]\times[-1,1]$ maps a sufficiently smooth function $u:[-1,1]^3 \to \R$ to
\begin{equation}\label{eq:LinearDifferentialOperator}
    \L u(x,y,z) = \sum_{a=0}^{N_x} \sum_{b=0}^{N_y} \sum_{c=0}^{N_z} \alpha_{abc}(x,y,z) \frac{\partial^{a+b+c}}{\partial x^a \partial y^b \partial z^c} u(x,y,z), 
\end{equation} 
where $N_x,N_y,N_z$ are called differential order and $\alpha_{abc}(x,y,z):[-1,1]^3 \to \R$ are coefficient functions for $\rev{0}\leq a \leq N_x,\ \rev{0} \leq b \leq N_y,\ \rev{0}\leq c \leq N_z$.
For the differential operator $\L$ we consider the linear PDE \begin{equation} \label{eq:PDEdefinition} \L u = f\end{equation} with right hand side $f:[-1,1]^3 \to \R$. 
The PDE can be solved uniquely when the system is complemented with sufficient boundary conditions. 

\subsection{Approximation format} 
We approximate the solution  $u:[-1,1]^3 \to \mathbb{R}$ of the PDE~\eqref{eq:PDEdefinition} in the space $\mathbb{P}_{n_1,n_2,n_3}$ of trivariate polynomials of degree at most $(n_1,n_2,n_3)$.
We express $u$ in terms of a \rev{tensorized} basis of Chebyshev polynomials, which leads to a representation of the form
\begin{equation}\label{eq:ChebyshevFormat}
    u(x,y,z) \approx  \sum_{i=\rev{0}}^{n_1} \sum_{j=\rev{0}}^{n_2} \sum_{k=\rev{0}}^{n_3} \mathcal{U}_{ijk} T_i(x) T_j(y) T_k(z),
\end{equation}
where $\mathcal{U}\in \R^{(n_1 \rev{+1)} \times (n_2 \rev{+1)} \times (n_3\rev{+1)}}$ is called \rev{the} coefficient tensor and $T_k(x) = \cos(\rev{k}\cos^{-1}(x))$ denotes the $k$-th Chebyshev polynomial \rev{for $x\in [-1,1]$}
\rev{Note that throughout this work, entries of vectors, matrices and tensors are indexed starting from $0$ following the notation in~\cite{NIST}.}
\section{Operator discretization}\label{sec:Discretization}
In the following, we discretize the differential operator $\L$ \rev{for fixed coefficient functions $\alpha_{abc}$} by approximating $\L$ as mapping from $\mathbb{P}_{n_1,n_2,n_3}$ to $\mathbb{P}_{n_1,n_2,n_3}$.
This is particularly easy for constant coefficients in the differential operator $\L$, i.e.\rev{,} $\alpha_{abc}(x,y,z) = \alpha_{abc} \in \R$ in Equation~\eqref{eq:LinearDifferentialOperator}. 
Applying a linear differential operator with constant coefficients to a polynomial does not increase the polynomial degree.
For non-constant coefficients an additional truncation is need\rev{ed} to obtain a polynomial in $\mathbb{P}_{n_1,n_2,n_3}$.
Thus, it is natural to see the application of the operator as a \rev{transformation} of the coefficient tensor.
We discuss the discretization for the constant  case first before generalizing to the non-constant case in Section~\ref{sec:NonConst}. 

\subsection{One-dimensional differential operators}\label{sec:1Doperators}
We briefly recapitulate how to obtain the mapping describing the \rev{transformation} of coefficients in a one-dimensional setting before returning to the three-dimensional setting.
Let $u : [-1,1] \to \R$ be a polynomial of degree $n$ represented in \rev{the} Chebyshev basis by 
$u(x) =  \sum_{k=\rev{0}}^n u_k T_k(x),$
with coefficients $u_k \in \R$ and Chebyshev polynomials $T_k(x)$.
Applying a one-dimensional linear differential operator $\L$ \rev{of order $N$} with constant coefficients  to $u$ can be written as 
\begin{equation} \label{eq:1Doperator}
\L u(x) = \sum_{a=0}^N \alpha_a \frac{d^a}{dx^a} u(x),
\end{equation}
with coefficients $\alpha_a \in \R$. This is a linear combination of (higher order) derivatives of $u$. 
For every derivative of the polynomial $u$, there exist \rev{a so-}called differentiation matrix which maps the coefficient vector $\mathbf{u} = (u_{\rev{0}},\dots,u_n)$ to the coefficients of the derivative. 
The remainder of this section follows the ideas of~\cite{Olver13} to represent the derivative using an ultraspherical basis instead of a Chebyshev basis, which leads to better conditioned and sparse differentiation matrices.

For the parameter $\lambda \rev{>0}$ and $k = \rev{0,1},\dots$, ultraspherical polynomials follow the recurrence relation
\[
\rev{
4 \lambda (k+\lambda+1) (1-x^2) C_k^{(\lambda+1)} (x) = -(k+1)(k+2) C_{k+2}^{(\lambda)}(x) + (k + 2 \lambda) (k + 2 \lambda + 1) C_k^{(\lambda)}(x),}
\]
\rev{where $C_k^{(1)}(x) = (\sin(k+1) \cos^{-1}(x)) / \sin(\cos^{-1}(x))$}~\cite{NIST}.
Let $v(x) = \sum_{k=\rev{0}}^n v_k C_k^{(\lambda)}(x)$ denote the $\lambda$th derivative of $u$ represented in $C^{(\lambda)}$ basis, then the coefficient vector $\mathbf{v} = (v_{\rev{0}},\dots,v_n)$ satisfies $\mathbf{v} = D_\lambda \mathbf{u}$, where the sparse differentiation matrix $D_\lambda\in \R^{(n\rev{+1)}\times (n\rev{+1)}}$ is defined as \[ D_\lambda = 2^{\lambda-1}(\lambda-1)! \begin{pmatrix}
\overbrace{0 \dots 0}^{\lambda \textrm{ times}} & \lambda & &  &\\
& & \lambda+1 &  &\\
& &  &  \ddots &\\
& &  & &  \rev{n} \\
& &  & &  \begin{array}{c}
0 \\ \vdots \\ 0
\end{array} 
\end{pmatrix}.\]
Note that the ultraspherical basis is different for different $\lambda$.
The sparse transformation matrices $S_0\in \R^{(n\rev{+1)}\times (n\rev{+1)}}$, mapping Chebyshev to $C^{(1)}$ coefficients, and $S_\lambda \in \R^{(n\rev{+1)}\times (n\rev{+1)}}$, mapping $C^{(\lambda)}$ to $C^{(\lambda+1)}$, are defined as
\[ 
{S}_0 = \begin{pmatrix}
1 & 0 & -\frac{1}{2} & & & &\\
& \frac{1}{2} & 0 & -\frac{1}{2} & & & \\
& & \frac{1}{2} & 0 & -\frac{1}{2}  & & \\
& &  & \ddots  & \ddots & \ddots &\\
& &  &  &  \rev{\frac{1}{2}} & \rev{0} & \rev{-\frac{1}{2}}\\
& &  &  & & \rev{\frac{1}{2}} & \rev{0} \\
& &  & & & & \rev{\frac{1}{2}}  \\
\end{pmatrix}, \quad
{S}_\lambda = \begin{pmatrix}
1 & 0 & -\frac{\lambda}{\lambda+2} & & & & \\
& \frac{\lambda}{\lambda+1} & 0 & -\frac{\lambda}{\lambda+3} & & & \\
& & \frac{\lambda}{\lambda+2} & 0 & -\frac{\lambda}{\lambda+4}  & &\\
& &  & \ddots  & \ddots & \ddots & \\
& &  &  &  \rev{\frac{\lambda}{\lambda+n-2}} & \rev{0} & \rev{-\frac{\lambda}{\lambda+n}}\\
& &  &  & & \rev{\frac{\lambda}{\lambda+n-1}} & \rev{0} \\
& &  & & & & \rev{\frac{\lambda}{\lambda+n}}  \\
\end{pmatrix}.
\]

Let $\L u(x) = \sum_{k=\rev{0}}^n w_k C^{(N)}_k(x)$ be represented in $C^{(N)}$ ultraspherical basis with coefficient vector $\mathbf{w} = (w_{\rev{0}},\dots,w_n)$, then
\begin{equation}\label{eq:DefinitionDiscreteL1D}
\mathbf{w} = (\underbrace{a_N D_N + a_{N-1} S_{N-1} D_{N-1} + \dots + a_1 S_{N-1} \cdots S_1 D_1  + a_0 S_{N-1} \cdots S_0}_{= L}) \mathbf{u}.
\end{equation}
The matrix $L \in \R^{(n\rev{+1)}\times (n\rev{+1)}}$ describes how the one-dimensional differential operator $\L$ acts on \rev{Chebyshev} coefficients \rev{of the solution}.

\subsection{Three-dimensional differential operators} \label{sec:3DSplittingConstantCoeffs}
In Chebop2~\cite{Townsend15} it is suggested to split two-dimensional differential operators via an SVD into a sum of one-dimensional operators, for which the matrices $L$  can be computed as in Section~\ref{sec:1Doperators}. Such splittings can be generalized to the three-dimensional setting using a CP decomposition. 

Let $\mathcal{A} \in \R^{(N_x+1)  \times (N_y+1) \times (N_z+1)}$ denote the tensor with entries $\mathcal{A}_{a,b,c}$ given by the coefficients $\alpha_{abc}$ in Equation~\eqref{eq:LinearDifferentialOperator} for $0\leq a \leq N_x,\ 0 \leq b \leq N_y,\ 0 \leq c \leq N_z$.
A CP decomposition~\cite{Kolda09} of $\mathcal{A}$ takes the form
\[
\mathcal{A} = \sum_{r=1}^R  \mathbf{a}^{(r)} \circ \mathbf{b}^{(r)} \circ \mathbf{c}^{(r)} \rev{,}
\]
where $R$ is called \rev{the} tensor-rank, $\circ$ denotes the outer product\rev{,} and  $\mathbf{a}^{(r)}\in \R^{N_x+1},\ \mathbf{b}^{(r)} \in \R^{N_y+1},\ \mathbf{c}^{(r)} \in \R^{N_z+1}$ for $r=1,\dots,R$.

We can now rewrite the application of $\mathcal{L}$ to $u$ as 
\begin{equation}\label{eq:DefSplitting1DOperators}
\mathcal{L}u(x) = \sum_{r=1}^R \underbrace{ \sum_{a = 0}^{N_x} \mathbf{a}^{(r)}_{a}  \frac{\partial^a}{\partial x^a} }_{= \mathcal{L}_r^{(x)}} 
\underbrace{ \sum_{b = 0}^{N_y} \mathbf{b}^{(r)}_{b}  \frac{\partial^b}{\partial y^b} }_{= \mathcal{L}_r^{(y)}}
\underbrace{ \sum_{a = 0}^{N_z} \mathbf{c}^{(r)}_{c}  \frac{\partial^c}{\partial z^c} }_{= \mathcal{L}_r^{(z)}}
u(x,y,z),
\end{equation}
in terms of one-dimensional differential operators $\L_r^{(x)},\L_r^{(y)},\L_r^{(z)}$.

Let  $L_r^{(x)}\in \R^{(n_1\rev{+1)} \times (n_1\rev{+1)}},\ L_r^{(y)} \in \R^{(n_2\rev{+1)} \times (n_2\rev{+1)}},\ L_r^{(z)} \in \R^{(n_3\rev{+1)} \times (n_3\rev{+1)}}$ denote the matrices $L$ associated with the operators $\L_r^{(x)},\L_r^{(y)},\L_r^{(z)}$ as defined in Equation~\eqref{eq:DefinitionDiscreteL1D}. Let $u$ be a polynomial with coefficient tensor $\mathcal{U} \in \R^{(n_1\rev{+1)} \times (n_2\rev{+1)} \times (n_3\rev{+1)}}$ satisfying Equation~\eqref{eq:ChebyshevFormat} with equality. Then 
\[
\mathcal L u(x,y,z) = \sum_{i=\rev{0}}^{n_1} \sum_{j=\rev{0}}^{n_2} \sum_{k=\rev{0}}^{n_3} \mathcal{V}_{ijk} C_i^{(N_x)}(x) C_j^{(N_y)} (y) C_j^{(N_z)}(z), 
\]
with coefficient tensor
\begin{equation}\label{eq:CoeffTensorEquation3D}
    \mathcal{V} = \sum_{r=1}^R \mathcal{U} \times_1 L_r^{(x)} \times_2 L_r^{(y)} \times_3 L_r^{(z)},
\end{equation}
where $\mathcal{V} \times_\ell M$ denotes the mode-$\ell$ multiplication. For a tensor $\rev{\mathcal{T}} \in \R^{n_1 \times n_2 \times n_3}$ and a matrix $M \in \R^{m \times n_{\ell}}$  \rev{the mode-$\ell$ multiplication $\times_{\ell}$} is defined as the multiplication of \rev{$M$ with} every mode-$\ell$ fiber of $\mathcal{T}$, i.e.
\begin{equation*}
    (\mathcal{T} \times_{\ell} M)^{\set{{\ell}}} = M \mathcal{T}^{\set{{\ell}}},
\end{equation*}
where $\mathcal{T}^{\set{{\ell}}}$ denotes the mode-${\ell}$ matricization, which is the matrix containing all mode-${\ell}$ fibers of $\mathcal{T}$ \rev{in its columns}~\cite{Kolda09}. 

\subsection{Generalization to non-constant coefficients}\label{sec:NonConst}
\subsubsection{One-dimensional differential operators}\label{sec:1dNonConst}
We now consider the one-dimensional differential operator~\eqref{eq:1Doperator} with non-constant coefficients $\alpha_a:[-1,1] \to \R$. In this section, we define multiplication matrices to incorporate these coefficients in the discretization.

Let $u(x)$ be a polynomial of the form 
$
u(x) = \sum_{k=\rev{0}}^n u_k B_k(x),
$
with basis functions $B_k(x)$ chosen as Chebyshev polynomials $T_k$ or ultraspherical polynomials $C_k^{(\lambda)}$. 
To multiply $u(x)$ \rev{by} a function $v:[-1,1]\to \R$, we approximate $v(x) \approx \tilde{v}(x) = \sum_{k=\rev{0}}^n v_k B_k(x)$ by a polynomial in the same basis using Chebyshev interpolation~\cite{Trefethen13} (and basis transformations). 
Further, we also approximate the product $\tilde{v}(x)u(x) \approx z(x) = \sum_{k=\rev{0}}^n z_k B_k(x)$ by a polynomial in the same basis. There exist so\rev{-}called multiplication matrices~\cite{Olver13,Fortunato21} depending on the coefficients $\mathbf{v} = (v_{\rev{0}},\dots,v_n)$, which map the coefficients $\mathbf{u} = (u_{\rev{0}},\dots,u_n)$ to the coefficients $\mathbf{z} = (z_{\rev{0}},\dots, z_n)$ such that $z(x)$ approximates $\tilde{v}(x)u(x)$. In the following, we summarize how these matrices are defined in~\cite{Townsend14} for both Chebyshev and ultraspherical bases. For the Chebyshev basis, we define $\mathbf{z} = M[\mathbf v] \mathbf{u}$ and the multiplication matrix $M[\mathbf{v}] \in \R^{(n\rev{+1)} \times (n\rev{+1)}}$ given by
\[
M[\mathbf{v}] =  \frac{1}{2} \begin{pmatrix} 2 v_{\rev{0}} & v_{\rev{1}} &v_{\rev{2}} & \dots & \rev{v_{n-1}} & \rev{v_n}\\
v_{\rev{1}}  & 2 v_{\rev{0}} & v_{\rev{1}}& \dots & \rev{v_{n-2}} & \rev{v_{n-1}}\\
v_{\rev{2}} & v_{\rev{1}} & 2 v_{\rev{0}} & \dots & \rev{v_{n-3}} & \rev{v_{n-2}}\\
\rev{\vdots} & \rev{\vdots} & \rev{\vdots} & \rev{\ddots} & \rev{\vdots}  & \rev{\vdots}\\
\rev{v_{n-1}} & \rev{v_{n-2}} & \rev{v_{n-3}} & \dots & \rev{2 v_0} & \rev{ v_1} \\
\rev{v_n} & \rev{v_{n-1}} & \rev{v_{n-2}} & \dots & \rev{v_1} & \rev{2 v_0}
\end{pmatrix} + \frac{1}{2}
\begin{pmatrix}
0 & 0 & 0 & \dots & \rev{0} & \rev{0} & \rev{0} \\
v_{\rev{1}} & v_{\rev{2}} & v_{\rev{3}} &  \dots & \rev{v_{n-1}} & \rev{v_n} & \rev{0} \\
v_{\rev{2}} & v_{\rev{3}} & v_{\rev{4}} & \dots & \rev{v_n} & \rev{0} & \rev{0}\\
\vdots & \vdots & \vdots & \ddots & \rev{\vdots} & \rev{\vdots} & \rev{\vdots} \\
\rev{v_{n-1}} & \rev{v_n}  & \rev{0} & \dots & \rev{0} & \rev{0} & \rev{0} \\
\rev{v_n} & \rev{0} & \rev{0} & \dots & \rev{0} & \rev{0} & \rev{0} \\
\end{pmatrix}.
\]
For the $C^{(\lambda)}$ ultraspherical basis, we define $\mathbf{z} = M^{(\lambda)}[\mathbf v] \mathbf{u}$ and the multiplication matrix $M^{(\lambda)}[\mathbf{v}] \in \R^{(n\rev{+1)} \times (n\rev{+1)}}$ as $M^{(\lambda)}[\mathbf{v}] = \sum_{i=\rev{0}}^n v_i M_{i}^{(\lambda)}$. 
The matrices $M_{i}^{\rev{(\lambda)}} \in \R^{(n\rev{+1)} \times (n\rev{+1)}}$ are defined recursively for $i = \rev{0,1},\dots$ by
 \[\rev{(i+2)}M_{i+2}^{(\lambda)} = 2(i+\lambda\rev{+1}) N^{(\lambda)} M^{(\lambda)}_{i+1} - (i+2\lambda) M^{(\lambda)}_{i},\]
with $M_{\rev{0}}^{(\lambda)} = I$, $M_{\rev{1}}^{(\lambda)} = 2 \lambda N^{(\lambda)}$ and
\begin{align*}
    N^{(\lambda)} = \begin{pmatrix}
    0 & \frac{2\lambda}{2(\lambda+1)} & & & & & \\
    \frac{1}{2\lambda} & 0 & \frac{2\lambda+1}{2(\lambda+2)} & & & &\\
    0 & \frac{2}{2(\lambda+1)} & 0 & \frac{2\lambda+2}{2(\lambda+3)} & & &\\
    & & \frac{2}{2(\lambda + 2)} & 0 & \rev{\frac{2\lambda+3}{2(\lambda+4)}} & & \\
    & &  & \ddots & \ddots & \rev{\ddots} & \\
    & & &  & \rev{ \frac{2}{2(\lambda + n-2)} }  & \rev{0} & \rev{ \frac{2\lambda+n-1}{2(\lambda+n)}} \\
    & & &  &  & \rev{ \frac{2}{2(\lambda + n-1)} }  & \rev{0} \\
    \end{pmatrix} \rev{\in \R^{(n+1) \times (n+1)}}.
\end{align*}

We now incorporate the non-constant coefficients $\alpha_a(x)$ in\rev{to} the coefficient mapping.
Let $\mathbf{a}_a \in \R^{(n\rev{+1)}}$ denote coefficients of polynomial approximations of $\alpha_a(x)$ in \rev{the} Chebyshev basis for $a=0$ and in $C^{(a)}$ ultraspherical basis for $a =1,2,\dots$. Analogous to Equation~\eqref{eq:DefinitionDiscreteL1D}, we approximate $\L u \approx \sum_{k=\rev{0}}^n w_k T_k(x)$ by computing the coefficients $\mathbf{w} = (w_{\rev{0}},\dots,w_n)$ defined as
\[
\mathbf{w} = (\underbrace{ M^{(N)}[\mathbf{a}_N] D_N + S_{N-1} M^{(N-1)}[\mathbf{a}_{N-1}] D_{N\rev{-1}} \rev{+} \dots +  S_{N-1} \cdots S_1 M^{(1)}[\mathbf{a}_1] D_1  +  S_{N-1} \cdots S_0 M[\mathbf{a}_0] }_{= L}) \mathbf{u}.
\]
 
\rev{\textit{Remark.} The interpolation of $v$ and the approximation of the product $\tilde{v}u$ introduce truncation errors. For sufficiently large $n$, these errors are close to machine precision. }

\subsubsection{Three-dimensional differential operators}\label{sec:Discretize3DNonConst}
We now consider the three-dimensional differential operator $\L$ as defined in Equation~\eqref{eq:LinearDifferentialOperator} with non-constant coefficients. 
We proceed by splitting this operator into one-dimensional operators with non-constant coefficients.  

Let the polynomial degree $(n_1,n_2,n_3)$ be chosen sufficiently large to accurately approximate the coefficient functions $\alpha_{abc}(x,y,z)$ using \rev{tensorized} Chebyshev interpolation in the form of Equation~\eqref{eq:ChebyshevFormat}. Let $\mathcal{B}^{(abc)} \in \R^{(n_1\rev{+1)} \times (n_2\rev{+1)}  \times (n_3\rev{+1)} }$ denote the coefficient tensors corresponding to $\alpha_{abc}(x,y,z)$ for each multi-index $abc$. We define the tensor $\mathcal{A} \in \R^{(N_x+1)\times  (n_1\rev{+1)}  \times (N_y+1) \times (n_2\rev{+1)}  \times (N_z+1) \times (n_3\rev{+1)} }$ with entries
\[
\mathcal{A}_{aibjck} = \mathcal{B}^{(abc)}_{ijk}
\]
for $0 \leq a \leq N_x,\ 0 \leq b \leq N_y,\ 0 \leq c \leq N_z,\ \rev 0 \leq i \leq n_1,\ \rev 0 \leq j \leq n_2,\ \rev 0 \leq k \leq n_3$. This approximates the differential operator in the sense that
\[
\mathcal{L} u(x,y,z) \approx \tilde{\mathcal{L}} u(x,y,z) =  \sum_{a=0}^{N_x} \sum_{i=\rev 0}^{n_1}\sum_{b=0}^{N_y} \sum_{j=\rev 0}^{n_2} \sum_{c=0}^{N_z} \sum_{k=\rev 0}^{n_3} \mathcal{A}_{ai bjck} T_i(x) T_j(y) T_k(z) \frac{\partial^{a+b+c}}{\partial x^a \partial y^b \partial z^c} u(x,y,z).
\]

We now reshape $\mathcal{A}$ into a tensor of order 3 in $\R^{(N_x+1)  (n_1\rev{+1)} \times (N_y+1)  (n_2\rev{+1)} \times (N_z+1)  (n_3\rev{+1)}}$ and compute a CP decomposition of the form
\[
\mathcal{A} = \sum_{r=1}^R \mathbf{a}^{(r)} \circ \mathbf{b}^{(r)} \circ \mathbf{c}^{(r)},
\]
where $\mathbf{a}^{(r)} \in \R^{(N_x+1) (n_1\rev{+1)}},\ \mathbf{b}^{(r)} \in \R^{(N_y+1) (n_2\rev{+1)}},\ \mathbf{c}^{(r)} \in \R^{(N_z+1) (n_3\rev{+1)}}$. We reshape the vectors $\mathbf{a}_r,\mathbf{b}_r,\mathbf{c}_r$ back into matrices in $\mathbf{A}^{(r)} \in \R^{(N_x+1) \times (n_1\rev{+1)}},\ \mathbf{B}^{(r)} \in \R^{(N_y+1) \times (n_2\rev{+1)}}$ and $\mathbf{C}^{(r)} \in \R^{(N_z+1)\times (n_3\rev{+1)}}$\rev{,} respectively.
Analogously to Equation~\eqref{eq:DefSplitting1DOperators}, we now define one-dimensional differential operators
\begin{align*}
\L_r^{(x)}u(x) = \sum_{a=0}^{N_x} \Big( \sum_{i=\rev 0}^{n_1} \mathbf{A}^{(r)}_{ai} T_i(x) \Big) \frac{\partial^a}{\partial x^a} u(x), \\
\L_r^{(y)}u(y) = \sum_{b=0}^{N_y} \Big( \sum_{j=\rev 0}^{n_2} \mathbf{B}^{(r)}_{bj} T_j(y) \Big) \frac{\partial^b}{\partial y^b} u(y), \\
\L_r^{(z)}u(z) = \sum_{c=0}^{N_z} \Big( \sum_{k=\rev 0}^{n_3} \mathbf{C}^{(r)}_{ck} T_k(z) \Big) \frac{\partial^c}{\partial z^c} u(z),
\end{align*}
which satisfy by construction
\begin{equation}
    \label{eq:3dRawDiscretization}
\tilde{\L}u(x,y,z) = \sum_{r=1}^R  (\L_r^{(x)} \rev{\otimes} \L_r^{(y)} \rev{\otimes} \L_r^{(z)}) u(x,y,z),
\end{equation}
\rev{where $\rev{\otimes}$ denotes the tensor product for linear operators.}
Note that the terms of the form $\sum_{i=\rev 0}^{n_1} \mathbf{A}^{(r)}_{ai} T_i(x)$ are univaritate functions. So, each of the one-dimensional differential operators fits into the setting of Section~\ref{sec:1dNonConst} and we can obtain matrices $L_r^{(x)}, L_r^{(y)}, L_r^{(z)}$ as in Equation~\eqref{eq:CoeffTensorEquation3D}.

\section{A spectral method for three-dimensional linear PDEs}\label{sec:SpectralMethod}
In this section, we present how to compute approximate solutions for PDEs of the form~\eqref{eq:PDEdefinition}. 
We again discretize the differential operator $\L$ as in Section~\ref{sec:Discretization}.
Additionally, we discretize the right hand side $f$ using a truncated expansion with ultraspherical basis functions of the form
 \[
 f(x,y,z) \approx \sum_{i=\rev 0}^{n_1} \sum_{j=\rev 0}^{n_2} \sum_{k=\rev 0}^{n_3} \mathcal{F}_{ijk} C_i^{(N_x)}(x) C_j^{(N_y)} (y) C_k^{(N_z)}(z), 
 \]
with coefficient tensor $\mathcal{F} \rev{\in \R^{(n_1+1)\times(n_2+1)\times(n_3+1)}}$.
The discretized PDE reads as tensor-valued linear system of the form
 \begin{equation}\label{eq:DiscretizedPDE}
     \sum_{r=1}^R \mathcal{U} \times_1 L_r^{(x)} \times_2 L_r^{(y)} \times_3 L_r^{(z)} = \mathcal{F}.
 \end{equation}
It remains to incorporate the boundary conditions.

\subsection{Boundary condition discretization}
Following the ideas of~\cite{Townsend15}, we can discretize commonly used boundary conditions for three-dimensional PDEs on cubes as constraints of the form 
\begin{equation}
    \label{eq:GeneralDiscretizedBoundaryConditions}
    \mathcal U \times_1 B_1 = \mathcal G_1, \quad \mathcal U \times_2 B_2 = \mathcal G_2, \quad \mathcal U \times_3 B_3 = \mathcal G_3,
\end{equation}
where the matrices $B_1 \in \R^{N_x \times (n_1\rev{+1)}}, B_2 \in \R^{N_y \times (n_2\rev{+1)}}, B_3 \in \R^{N_z \times (n_3\rev{+1})}$ have linearly independent rows and $\rev{\mathcal G}_1 \in \R^{N_x \times (n_2\rev{+1)} \times (n_3\rev{+1)}},\rev{\mathcal G}_2 \in \R^{(n_1\rev{+1)} \times N_y \times (n_3\rev{+1)}},\rev{\mathcal G}_3 \in \R^{(n_1\rev{+1)} \times (n_2\rev{+1)} \times N_{\rev{z}}}$. 
We present two examples \rev{of} how constraints of the form~\eqref{eq:GeneralDiscretizedBoundaryConditions} can be derived. 
\rev{We want to emphasize that the discretized boundary conditions need to satisfy compatibility constraints as in~\cite[Section 6.4]{Townsend14}.}

\subsubsection{Dirichlet conditions}\label{sec:Dirichlet}
We consider the Dirichlet boundary condition $u(1,y,z) = h(y,z)$ for a given function $h:[-1,1] \to \R$. We approximate the function $h$ using bivariate Chebyshev interpolation in $(n_2\rev{+1)} \times (n_3\rev{+1)}$ points. Let $H \in \R^{(n_2\rev{+1)} \times (n_3\rev{+1)}}$ denote the corresponding coefficient \rev{matrix of $h$}. We can enforce that the solution $u$ given in Chebyshev basis~\eqref{eq:ChebyshevFormat} coincides with the interpolant of $h$ by demanding that \begin{equation}
\label{eq:DirichletConstruction}    
u(1,y,z) = \sum_{i=\rev 0}^{n_1} \sum_{j=\rev 0}^{n_2} \sum_{k=\rev 0}^{n_3} \mathcal{U}_{ijk} T_i(1) T_j(y) T_k(z) = \sum_{j=\rev 0}^{n_2} \sum_{k=\rev 0}^{n_3} H_{jk} T_j(y) T_k(z).
\end{equation}
This can be equivalently written as $U \times_1 B_1 = \rev{\mathcal G}_1$ with $B_{\rev{1}} = (T_{\rev{0}}(1), \dots, T_{n_1}(1))$ and $\mathcal (G_1)_{\rev{ 0,0 : n_2, 0:n_3}} = H$.

We now consider an additional Dirichlet boundary condition $u(-1,y,z) = \tilde{h}$. Again, we compute the coefficient \rev{matrix} $\tilde{H}$ corresponding to $\tilde{h}$.  Analogously to Equation~\eqref{eq:DirichletConstruction}, we enforce that $u(-1,y,z)$ coincides with the interpolant of $\tilde{h}$. We can express both conditions simultaneously in the form of~\eqref{eq:GeneralDiscretizedBoundaryConditions} by defining
\[ 
 B_1 = \begin{pmatrix} T_{\rev{0}}(-1) & T_{\rev{1}}(-1) & \dots & T_n(-1) \\ T_{\rev{0}}(1) & T_{\rev{1}}(1) & \dots & T_n(1) \end{pmatrix},
\]
and $(G_1)_{\rev{ 0,0 : n_2, 0:n_3}} = \tilde{H},\ (G_{\rev 1})_{\rev{ 1,0 : n_2, 0:n_3}} = H$.

\subsubsection{\rev{Mixed Dirchlet and} Neumann conditions}\label{sec:Neumann}
We consider \rev{ one example of mixed boundary conditions with} Neumann boundary conditions \rev{on the right side of the cube $[-1,1]^3$ and Dirichlet boundary conditions on all other sides. The Neumann boundary condition is given by} $\frac{\partial}{\partial x} u(1,y,z) = h(y,z)$ for a given function $h$ and $u \in \mathbb{P}_{n_1\rev{,}n_2\rev{,}n_3}$ represented in Chebyshev basis~\eqref{eq:ChebyshevFormat}. As in the Dirichlet case, we use bivariate Chebyshev interpolation $h$ to obtain the coefficient \rev{matrix} $H \in \R^{(n_2\rev{+1)} \times (n_3\rev{+1)}}$. We now demand that 
\begin{equation}
    \label{eq:NeumannBCS}
    \frac{\partial}{\partial x} u(1,y,z) = \sum_{i=\rev 0}^{n_1} \sum_{j=\rev{0}}^{n_2} \sum_{k=\rev 0}^{n_3} \mathcal{U}_{ijk} \rev{T_i'}(1) T_j(y) T_k(z) = \sum_{j=\rev 0}^{n_2} \sum_{k=\rev 0}^{n_3} H_{jk} T_j(y) T_k(z),
\end{equation} 
where \rev{$T_i'(1) = i^2$ for $i=0,\dots,n$~\cite{Trefethen13}.} 
Equation~\eqref{eq:NeumannBCS} can be expressed in the form of~\eqref{eq:GeneralDiscretizedBoundaryConditions} with \[B_1 = (\rev{T_{\rev{0}}'}(1),\dots,\rev{T_{n_1}'}(1)) \quad \text{ and } \quad \mathcal (G_1)_{\rev{0},0:n_2,0:n_3} = H.\]
\rev{The Dirichlet boundary conditions on the left side can be included in $B_1,\mathcal{G}_1$ as in Section~\ref{sec:Dirichlet}.}

\subsection{Incorporating the boundary conditions}
We need to incorporate the discretized boundary conditions~\eqref{eq:GeneralDiscretizedBoundaryConditions} into the discretized PDE~\eqref{eq:DiscretizedPDE} to  obtain the unique solution $\mathcal{U}$. 
Following the ideas of~\cite[Section~6]{Townsend15}, we compute $\mathcal{U}$ by substituting~\eqref{eq:GeneralDiscretizedBoundaryConditions} into~\eqref{eq:DiscretizedPDE}.

Since $B_1,B_2,B_3$ have linearly independent rows, we can assume without loss of generality that 
\begin{equation}
    \label{eq:BoundaryEquationFormAssumption}
    (B_1)_{\rev{0} : (N_x \rev{-1)},\rev{0} : (N_x \rev{-1)}} = I,\quad (B_2)_{\rev{0} :(N_y \rev{-1)},\rev{0} : (N_y \rev{-1)}} = I,\quad (B_3)_{(\rev{0}:(N_z \rev{-1)},\rev{0} : (N_z \rev{-1)}} = I.
\end{equation}

\rev{For $r=1,\dots,R$, we rewrite the boundary conditions as}
\begin{align*}
\rev{-\mathcal U \times_1 (L_r^{(x)})_{\rev{0 : n_1,\rev{0} : (N_x \rev{-1)}}} B_1} &\rev{= -\mathcal G_1 \times_1 (L_r^{(x)})_{\rev{0 : n_1,\rev{0} : (N_x \rev{-1)}}}}, \\
\rev{-\mathcal U \times_2 (L_r^{(y)})_{\rev{0 : n_2,\rev{0} : (N_y \rev{-1)}}} B_2} &\rev{= -\mathcal G_2 \times_2 (L_r^{(y)})_{\rev{0 : n_2,\rev{0} : (N_y \rev{-1)}}}}, \\
\rev{-\mathcal U \times_3 (L_r^{(z)})_{\rev{0 : n_3,\rev{0} : (N_z \rev{-1)}}} B_3} &\rev{= -\mathcal G_3 \times_3 (L_r^{(z)})_{\rev{0 : n_3,\rev{0} : (N_z \rev{-1)}}}}.
\end{align*}
Substituting these \rev{modified} boundary conditions \rev{into the discretized PDE~\eqref{eq:DiscretizedPDE}} leads to
\begin{align*}
    &\sum_{r=1}^R \mathcal{U} \times_1 \underbrace{(L_r^{(x)}-(L_r^{(x)})_{\rev{0 : n_1,\rev{0} : (N_x \rev{-1)}}}B_1)}_{= \tilde {L}_r^{(x)} } \times_2 \underbrace{(L_r^{(y)}-(L_r^{(y)})_{\rev{ 0 : n_2,\rev{0} : (N_y \rev{-1)}}}B_2)}_{= \tilde{L}_r^{(y)}} \times_3 \underbrace{(L_r^{(z)}-(L_r^{(z)})_{\rev{0:n_3,\rev{0} : (N_z \rev{-1)}}}B_3)}_{= \tilde{L}_r^{(z)}} = \\
&\quad \mathcal{F} - \sum_{r=1}^R \rev{\mathcal G}_1 \times_1 (L_r^{(x)})_{\rev{(0:n_1, \rev{0} : (N_x \rev{-1)}}} \times_2 L_r^{(y)} \times_3 L_r^{(z)} \\ &\quad - \sum_{r=1}^R \rev{\mathcal G}_2 \times_1 (L_r^{(x)}-(L_r^{(x)})_{\rev{0:n_1, \rev{0} : (N_x \rev{-1)}}}B_1) \times_2 (L_r^{(y)})_{\rev{0:n_2,\rev{0} \colon (N_y \rev{-1)}}} \times_3 L_r^{(z)}  \\ &\quad - \sum_{r=1}^R \rev{\mathcal G}_3 \times_1 (L_r^{(x)}-(L_r^{(x)})_{\rev{0:n_1,\rev{0} : (N_x \rev{-1)}}}B_1) \times_2 (L_r^{(y)}-(L_r^{(y)})_{\rev{0:n_2 ,\rev{0} : (N_y \rev{-1)}}}B_2) \times_3 (L_r^{(z)})_{\rev{0:n_3,\rev{0} : (N_z \rev{-1)}}},
\end{align*}
where we denote the right hand side by $\tilde{\mathcal F}$. Observe that the first $N_x,N_y,N_z$ columns of $\tilde{L}_r^{(x)},\tilde{L}_r^{(y)},\tilde{L}_r^{(z)}$ are zero due assumption~\eqref{eq:BoundaryEquationFormAssumption}.
\begin{figure}
    \centering
    \includegraphics{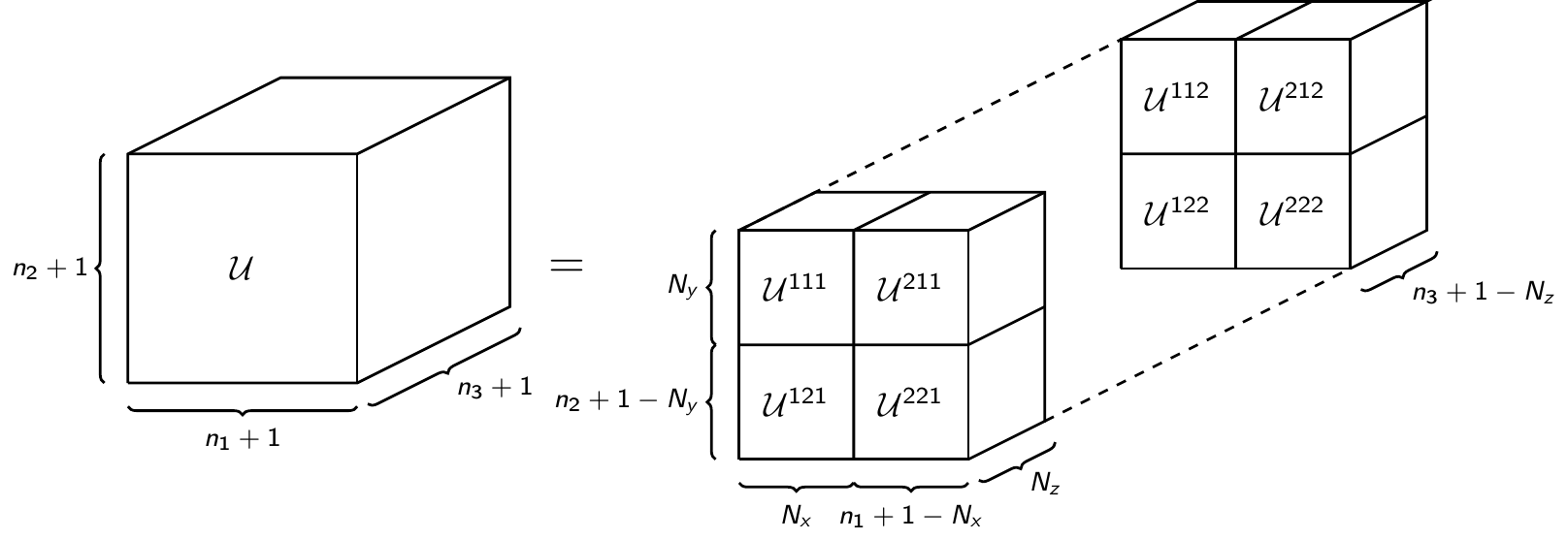}
    \caption{Visualization of the block decomposition of $\mathcal U \in \R^{(n_1\rev{+1)} \times (n_2\rev{+1)} \times (n_3\rev{+1)}}$ with block $\rev{\mathcal{U}^{111}} \in \R^{N_x \times N_y \times N_z}$.}
    \label{fig:BlockDecomposition}
\end{figure}
Let $\mathcal{U}$ be decomposed into tensor blocks as described in Figure~\ref{fig:BlockDecomposition}. The system after substitution uniquely determines the block $\rev{\mathcal{U}^{222}} = \rev{\mathcal{U}_{N_x: n_1,N_y: n_2,N_z: n_3}}$, which can be written as 
\begin{align}
     &\sum_{r=1}^R \rev{\mathcal{U}^{222}} \times_1 \underbrace{(\tilde{L}_r^{(x)})_{\rev{ 0: (n_1-N_x),N_x: n_1}}}_{= \hat{L}_r^{(x)}} \times_2 \underbrace{(\tilde{L}_r^{(y)})_{\rev{0 : (n_2-N_y),N_y : n_2}}}_{= \hat{L}_r^{(y)}} \times_3   \underbrace{(\tilde{L}_r^{(z)})_{\rev{ 0 : (n_3-N_z),N_z : n_3}}}_{=\hat{L}_r^{(z)}} \nonumber \\ & \quad  =  \underbrace{\tilde{\mathcal{F}}_{\rev{ \rev{0} : (n_1-N_x),\rev 0 : (n_2-N_y),\rev 0 : (n_3-N_z)}}}_{=\hat{\mathcal F}}. \label{eq:3DReducedSystem}
\end{align}
The remaining blocks can be reconstructed from the discretized boundary conditions~\eqref{eq:GeneralDiscretizedBoundaryConditions} as
\rev{\begin{align*}
    \mathcal{U}^{122} &= \rev{\mathcal G}_1 - \mathcal{U}^{222} \times_1 B_1, & \mathcal{U}^{212} &= \rev{\mathcal G}_2 - \mathcal{U}^{222} \times_2 B_2,\\
\mathcal{U}^{221} &= \rev{\mathcal G}_3 - \mathcal{U}^{222} \times_3 B_3, & \mathcal{U}^{112} &= \rev{\mathcal G}_1 - \mathcal{U}^{212} \times_1 B_1,\\
\mathcal{U}^{121} &= \rev{\mathcal G}_1 - \mathcal{U}^{221} \times_1 B_1, & \mathcal{U}^{211} &= \rev{\mathcal G}_3 - \mathcal{U}^{212} \times_3 B_3,  \\
\mathcal{U}^{111} &= \rev{\mathcal G}_1 - \mathcal{U}^{211} \times_1 B_1. & & 
\end{align*}}

\subsection{Solving tensor-valued linear systems}

The computation of $\mathcal{U}$ requires the solution of the unconstrained tensor-valued linear system~\eqref{eq:3DReducedSystem} of the form
\begin{equation}\label{eq:Structure3DReducedSystem}
    \sum_{r=1}^R \rev{\mathcal{U}^{222}} \times_1 \hat{L}_r^{(x)} \times_2 \hat{L}_r^{(y)} \times_3 \hat{L}_r^{(z)} = \hat{\mathcal F}.
\end{equation}
This system can be solved by reshaping the tensor-valued linear system into a vector-valued linear system of the form 
\begin{equation}
    \label{eq:ReshapedSystem}
\Big(\sum_{r=1}^R \hat{L}_r^{(z)} \otimes \hat{L}_r^{(y)} \otimes \hat{L}_r^{(x)}\Big) \textrm{vec}(\rev{\mathcal{U}^{222}}) = \textrm{vec}(\hat{\mathcal{F}}),
\end{equation}
where $\otimes$ denotes the Kronecker product and $\textrm{vec}(\cdot)$ denotes the vectorization~\cite{Kolda09}. 
For certain PDEs we can transform the system~\eqref{eq:Structure3DReducedSystem} into a Laplace-like equation
\begin{equation} \label{eq:LaplaceLikeEquation}
\rev{\mathcal{U}^{222}} \times_1 U  + \rev{\mathcal{U}^{222}}  \times _2  V + \rev{\mathcal{U}^{222}} \times_3  W =  \breve{\mathcal{F}},    
\end{equation}
with matrices $U \in \R^{(n_1\rev{+1}-N_x) \times (n_1\rev{+1}-N_x)},V \in \R^{(n_2\rev{+1}-N_y) \times (n_2\rev{+1}-N_y)},W \in \R^{(n_3\rev{+1}-N_z) \times (n_3\rev{+1}-N_z)}$ and tensor $\breve{\mathcal{F}} \in \R^{(n_1 \rev{+1}- N_x) \times (n_2\rev{+1}-N_y) \times (n_3 \rev{+1} - N_z)}$.
For instance, this can be achieved for $N_x=N_y=N_z$, $n_1=n_2=n_3$ and $B_1= B_2 = B_3$, when there exists a CP decomposition of the tensor $\mathcal{A}$ defined in Section~\ref{sec:Discretization} with symmetry constraints~\cite{Carroll80} of the form
\begin{equation}\label{eq:SpecialCPformForLaplacelike}
\mathcal{A} = \mathbf{a} \circ \mathbf{v} \circ \mathbf{w} + \mathbf{u} \circ \mathbf{b} \circ \mathbf{w} + \mathbf{u} \circ \mathbf{v} \circ \mathbf{c}.   
\end{equation}
Then Equation~\eqref{eq:Structure3DReducedSystem} is equivalent to the Laplace-like equation~\eqref{eq:LaplaceLikeEquation} with 
\begin{align}
    \label{eq:DefinitionMatricesLaplaceLike}
    U &= (\hat{{L}}_1^{(y)} )^{-1} \hat{{L}}_1^{(x)},\quad V = (\hat{{L}}_2^{(x)} )^{-1} \hat{{L}}_2^{(y)},\quad W = (\hat{{L}}_3^{(x)} )^{-1} \hat{{L}}_3^{(z)}, \nonumber \\ \breve{\mathcal{F}}  &= \hat{\mathcal{F}} \times_1 (\hat{{L}}_1^{(y)} )^{-1} \times_2 (\hat{{L}}_2^{(x)} )^{-1} \times_3 (\hat{{L}}_3^{(x)} )^{-1}.
\end{align}

To solve Laplace-like equations~\eqref{eq:LaplaceLikeEquation}, we apply the recursive blocked algorithm developed in~\cite{Chen19}. It transforms the matrices $U,V,W$ into quasi-triangular form by computing Schur decompositions. Block decompositions for the quasi-triangular matrices reveal an equivalent system of  Laplace-like equations with smaller matrices.
\rev{We apply this observation recursively}, until we can solve the small Laplace-like equations efficiently by reshaping\rev{. This} yields the blocked recursive algorithm. For $n = n_1=n_2=n_3$ this approach has a theoretical runtime of $\mathcal{O}(n^4)$ operations. In Section~\ref{sec:RuntimeComparison}, we demonstrate that this recursive blocked algorithm is much faster than directly reshaping the tensor-valued linear system.

\textit{Remark.} \rev{
For general PDEs with the same type of boundary conditions in each mode, we can obtain a CP decomposition of the form~\eqref{eq:SpecialCPformForLaplacelike} if the differential operator does not contain mixed derivatives and if all coefficients $\alpha_{abc}$ only depend on the variable corresponding to the mode in which their corresponding derivative acts. This includes, for instance, differential operators of the form
\[
\alpha_{1}(x) + \alpha_{2}(x) \frac{\partial}{\partial x}+ \alpha_{3}(x) \frac{\partial^2}{\partial x^2} + \alpha_{4}(y) \frac{\partial}{\partial y}+ \alpha_{5}(y) \frac{\partial^2}{\partial y^2} + \alpha_{6}(z) \frac{\partial}{\partial z}+ \alpha_{7}(z) \frac{\partial^2}{\partial z^2},
\]
with univariate coefficient functions  $\alpha_i:[-1,1] \to \R$.
}

\textit{Remark.} \rev{So far, we introduced the global spectral method for fixed polynomial degrees $(n_1,n_2,n_3)$. In order to heuristically determine if the solution is accurate, we can analyze the residual of~\eqref{eq:DiscretizedPDE} and the decay of the coefficients in $\mathcal{U}$~\cite{Aurentz17}. This can be used to adaptively increase $(n_1,n_2,n_3)$ until the solution is accurate.}

\textit{Remark.} \rev{In general, it is not recommended to explicitly invert the matrices in~\eqref{eq:DefinitionMatricesLaplaceLike}. The inversion can be avoided completely by extending the recursive solver in~\cite{Chen19} similar to how recursive solvers for Sylvester equations can be extended to solve generalized Sylvester equations.}
\section{Numerical results}\label{sec:NumericalExperiments}
All numerical experiments in this section were performed in MATLAB R2018b on a Lenovo Thinkpad T480s with Intel Core i7-8650U CPU and 15.4 GiB RAM. The code to reproduce these results is available from \url{https://github.com/cstroessner/SpectralMethod3D}.

\subsection{Runtime comparison} \label{sec:RuntimeComparison}
We consider Poissons's equation $\Delta u = f$ with homogeneous \rev{zero} Dirichlet boundary conditions.
In Table~\ref{tab:RuntimeComparison}, we compare our global spectral method to the nested alternating direction implicit method~(NADIM) proposed in~\cite[Section~5]{Fortunato19}. \rev{NADIM relies on solving two-dimensional Sylvester equations recursively on three levels. On each level an iterative algorithm is used, which leads to a large total number of iterations.}

\begin{table}[!ht]
\begin{tabular}{l|l|l|l|l|l|l|l|l|l}
          & \multicolumn{2}{c|}{$n = 10$} & \multicolumn{2}{c|}{$n = 30$} & \multicolumn{2}{c|}{$n = 50$}         & \multicolumn{2}{c|}{$n = 150$} & \multicolumn{1}{c}{} \\ \cline{2-9}
          & Time        &  Error        &   Time        &   Error        &   Time &   Error &   Time         &   Error        &                      \\ \cline{1-9} 
NADIM &   $21.8$        &  $1.47 \cdot 10^{-5}$ &   $742$      & $2.80 \cdot 10^{-8}$      &  -               &  -            &         -     &  - &          \\
reshape     &     $0.019$      & $1.55 \cdot 10^{-5}$  &  $2.40$       & $1.22 \cdot 10^{-15}$      &     $102$           &   $1.22 \cdot 10^{-15}$   & - &    -    &   \\
recursive   &  $0.054$         &  $1.55 \cdot 10^{-5}$      & $0.10$        & $3.12 \cdot 10^{-13}$      & $0.41$ & $8.20 \cdot 10^{-13}$ &     $10.9$            &  $1.15 \cdot 10^{-10}$             &           
\end{tabular}
\caption{Comparison of algorithms to solve Poisson's equation with known solution $u^*(x,y,z) = \sin(\pi x) \sin(\pi y) \sin(\pi z)$, \rev{from which the right hand side $f$ is explicitly computed.} To compute solutions $u \in \mathbb{P}_{n,n,n}$ with our global spectral method\rev{,} we compare reshaping~\eqref{eq:3DReducedSystem} and solving~\eqref{eq:ReshapedSystem} with backslash in \rev{MATLAB} (reshape) to transforming~\eqref{eq:3DReducedSystem} into a Laplace-like equation~\eqref{eq:LaplaceLikeEquation} and solving with the blocked recursive solver~\cite{Chen19} (recursive). Additionally, we compare to NADIM~\cite{Fortunato19}. For various $n$, we measure the runtime in seconds and we estimate  $||u-u^*||_\infty$, where $||\cdot||_\infty$ denotes the uniform norm, by sampling $1\,000$ random points.}
\label{tab:RuntimeComparison}
\end{table}

We observe that even though  NADIM has an asymptotic runtime of $\mathcal{O}(n^3 (\log(n))^3)$~\cite{Fortunato19}, it is the slowest algorithm in our setting and can not handle $n>30$ in a reasonable amount of time.  The asymptotic runtime of the recursive algorithm is slower with $\mathcal{O}(n^4)$, but in our experiments it is the fastest method and it can handle values of up to $n=150$ in less than $11$ seconds. For $n=50$ solving the reshaped system~\eqref{eq:ReshapedSystem} with backslash leads to an error of order $10^{-15}$, whereas the recursive algorithm only achieves an error of order $10^{-13}$.  While the recursive \rev{approach} is able to solve much larger systems, it is slightly more sensitive to numerical rounding errors.

\textit{Remark.} \rev{We want to emphasize that our approach is based on computing the coefficient tensor $\mathcal{U}$ fully. 
For certain PDEs including Poisson's equation, $\mathcal{U}$ potentially admits accurate low-rank approximations.
Computing a low-rank approximation directly as in~\cite{Shi21} can be much faster than computing the full coefficient tensor.}

\subsection{Stationary problems}

\subsubsection{Helmholtz problems}\label{sec:Helmholtz}
We consider the Helmholtz equation \rev{on $[-1,1]^3$ with non-homogeneous Dirichlet boundary conditions} as in~\cite{Hao16,Townsend15}, which arises for instance in the context of three-dimensional wave equations in acoustics~\cite{Wang97} and seismic-imaging~\cite{Plessix07}. It is defined as
\[
\Delta u(x,y,z) + \kappa^2  u(x,y,z) = f(x,y,z),
\] 
\rev{with} coefficient $\kappa \in \R$.  
For this differential operator, we define the tensor $\mathcal{A}$ as in Section~\ref{sec:3DSplittingConstantCoeffs}.  We observe that a CP decomposition in the form~\eqref{eq:SpecialCPformForLaplacelike} is given by
\begin{equation} \label{eq:SplittingHelmholtz}
\mathcal{A} =  
\begin{pmatrix} \kappa^2 \\ 0 \\ 1 \end{pmatrix} \circ
\begin{pmatrix} 1 \\ 0 \\ 0 \end{pmatrix} \circ
\begin{pmatrix} 1 \\ 0 \\ 0 \end{pmatrix} +
\begin{pmatrix} 1 \\ 0 \\ 0 \end{pmatrix} \circ
\begin{pmatrix} 0 \\ 0 \\ 1 \end{pmatrix} \circ
\begin{pmatrix} 1 \\ 0 \\ 0 \end{pmatrix} +
\begin{pmatrix} 1 \\ 0 \\ 0 \end{pmatrix} \circ
\begin{pmatrix} 1 \\ 0 \\ 0 \end{pmatrix} \circ
\begin{pmatrix} 0 \\ 0 \\ 1 \end{pmatrix}.
\end{equation}
Hence, we can derive a Laplace-like structure~\eqref{eq:LaplaceLikeEquation} for Equation~\eqref{eq:3DReducedSystem} and employ the recursive solver. 

In the following, we employ the global spectral method for the Helmholtz equation in~\cite[Section 5.3]{Canino98} given by 
\begin{equation}\label{eq:HelmholtzEquati}
    \Delta u(x,y,z) + \kappa(x)^2 u(x,y,z) = f(x,y,z),
\end{equation}
with function $\kappa(x) = \rev{\gamma_1} - \rev{\gamma_2} \cos(\pi \rev{\gamma_3} x / 2)$ and \rev{scalar} coefficients $\rev{\gamma_1,\gamma_2,\gamma_3} \in \R$. The right hand side $f$ and the  Dirichlet boundary conditions are \rev{computed explicitly from} the solution  
\begin{equation}\label{eq:HelmholtzSolution}
    u^*(x,y,z) = \exp( - \kappa(x)/\rev{\gamma_3} ) \cos( \pi \rev{\gamma_1} y/2) \cos(\pi \rev{\gamma_2} z/2).
\end{equation}
In order to incorporate the coefficient function $\kappa(x)$, we compute one-dimensional differential operators as in Equation~\eqref{eq:DefSplitting1DOperators} from the CP-decomposition~\eqref{eq:SplittingHelmholtz}. We then set  $\mathcal{L}^{(x)}_1 u(x) = \kappa(x)^2 u(x) + \frac{\partial^2}{\partial x^2}u(x)$ and discretize this operator as described in Section~\ref{sec:1dNonConst}. The resulting equation~\eqref{eq:3DReducedSystem} can still be transformed into a Laplace-like equation~\eqref{eq:LaplaceLikeEquation} by setting $U,V,W,\Breve{\mathcal{F}}$ as defined in~\eqref{eq:DefinitionMatricesLaplaceLike}.

\begin{figure}[ht]
    \centering
\subfloat{{\includegraphics[width=0.44\textwidth]{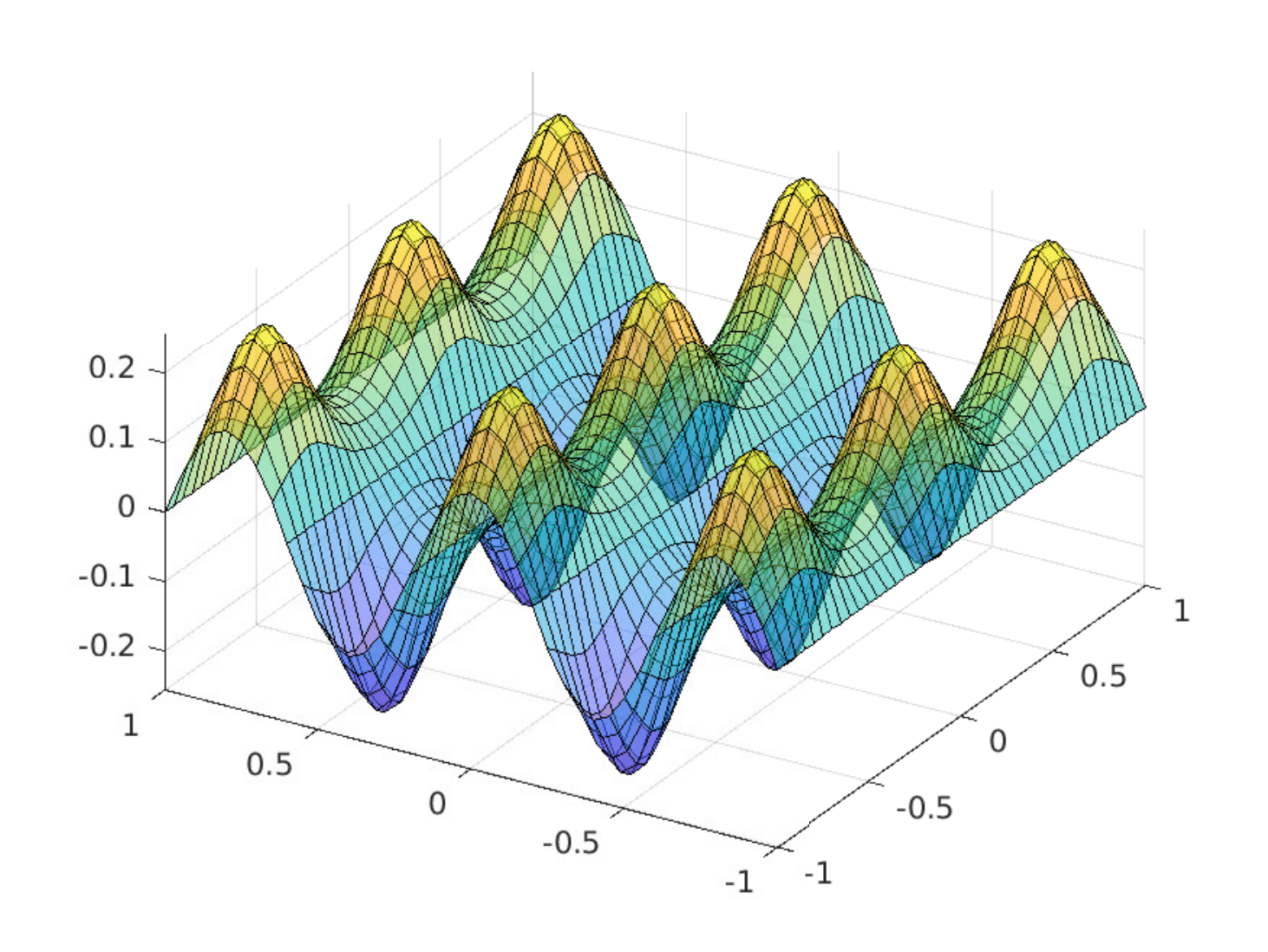}}}%
\qquad
\subfloat{{\includegraphics[width=0.44\textwidth]{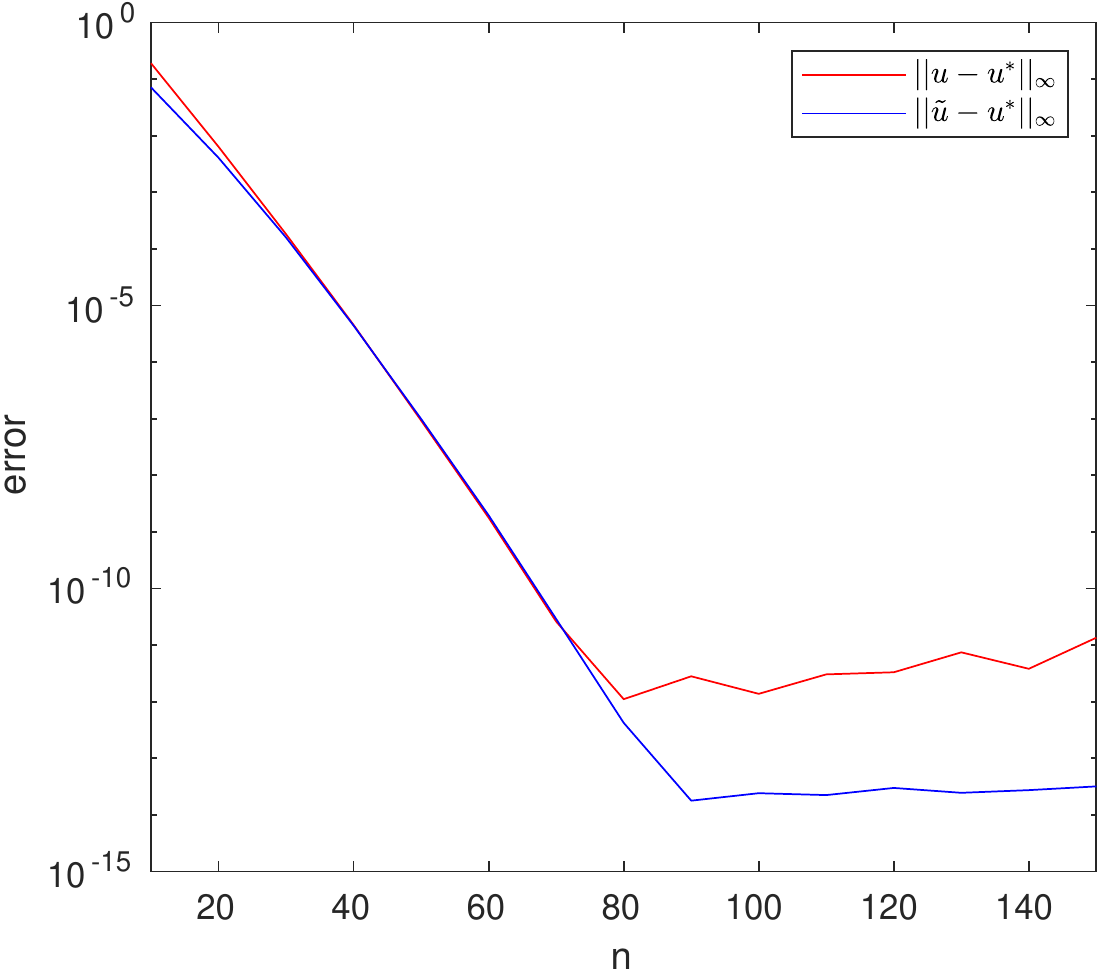}}}%
    \caption{ Left: We plot $u^*(x,y,z)$ as defined in~\eqref{eq:HelmholtzSolution} for $(\rev{\gamma_1,\gamma_2,\gamma_3}) = (5,3,5)$ and fixed $z=1/4$. Right: 
    We solve the Helmholtz equation~\eqref{eq:HelmholtzEquati} for $(\rev{\gamma_1,\gamma_2,\gamma_3}) = (5,3,5)$ with known solution $u^*$~\eqref{eq:HelmholtzSolution} via our global spectral method to obtain the solution $u \in \mathbb{P}_{n,n,n}$ for various $n$. Further, we compute the trivariate Chebyshev interpolant $\tilde{u} \in \mathbb{P}_{n,n,n}$ of $u^*$.
    For each $n$, we estimate $||u-u^*||_\infty$ and $||\tilde{u}-u^*||_\infty$ using $1\,000$ sample points (evaluation error). }
    \label{fig:Helmholtz}
\end{figure}

In Figure~\ref{fig:Helmholtz}, we depict how well the solution $u \in \mathbb{P}_{n,n,n}$ of our global spectral method approximates $u^*$.
We compare to the trivariate Chebyshev interpolant $\tilde{u}$ of $u^*$, which is close to the best approximation of $u$ in $\mathbb{P}_{n,n,n}$ and converges quickly~\cite{Stroessner21,Trefethen13}.
For $n<80$ the solution $u$ and the interpolant $\tilde{u}$ almost coincide. But for $n>80$ the error $||u-u^*||_\infty$ stagnates, whilst the interpolation error $||\tilde u - u^*||_\infty$ continues to decrease further before reaching a plateau close to machine precision. This discrepancy is \rev{again} caused by the \rev{the recursive blocked solver}.

\textit{Remark.} The derivation of the Laplace-like system from the CP decomposition can be extended from the Helmholtz equation to convection diffusion problems of the form
\[ - \nu \bigtriangleup u + \mathbf{\rev{\xi}}^T \bigtriangledown u = f, \]
with $\nu \in R$ and $\mathbf{\rev{\xi}} \in \R^3$ as studied in~\cite{Ballani13b, Chen12, Xiang10b}. For these problems we define the CP-decomposition 
\begin{equation*}
\mathcal{A} =  
\begin{pmatrix} 0 \\ \mathbf{\rev \xi}_1 \\ -\nu \end{pmatrix} \circ
\begin{pmatrix} 1 \\ 0 \\ 0 \end{pmatrix} \circ
\begin{pmatrix} 1 \\ 0 \\ 0 \end{pmatrix} +
\begin{pmatrix} 1 \\ 0 \\ 0 \end{pmatrix} \circ
\begin{pmatrix} 0 \\ \mathbf{\rev \xi}_2 \\ -\nu \end{pmatrix} \circ
\begin{pmatrix} 1 \\ 0 \\ 0 \end{pmatrix} +
\begin{pmatrix} 1 \\ 0 \\ 0 \end{pmatrix} \circ
\begin{pmatrix} 1 \\ 0 \\ 0 \end{pmatrix} \circ
\begin{pmatrix} 0 \\ \mathbf{ \rev \xi}_3 \\ -\nu \end{pmatrix}.
\end{equation*}

\subsubsection{Diffusion problems with separable coefficient}\label{sec:SeparableDiffusion}
Diffusion problems of the form 
\begin{equation}\label{eq:DiffusionProblem}
- \bigtriangledown \cdot (a(x,y,z) \bigtriangledown u(x,y,z)) = f(x,y,z),
\end{equation} 
with a separable coefficient $a(x,y,z) = a_1(x)a_2(y)a_3(z)$ defined by univariate functions $a_1,a_2,a_3: [-1,1] \to \R$, are not directly given as linear partial differential operator of the form~\eqref{eq:LinearDifferentialOperator}. We can, however, decompose the three-dimensional differential operator into a sum of three differential operators similar to~\eqref{eq:DefSplitting1DOperators} as
\[
\bigtriangledown \cdot (a(x,y,z) \bigtriangledown u(x,y,z)) =\Big( 
\underbrace{\frac{\partial}{\partial x}  \big(a(x,y,z) \frac{\partial}{\partial x} \big)}_{= \tilde{\mathcal{L}}_1^{(x)}} +
\underbrace{\frac{\partial}{\partial y} \big( a(x,y,z) \frac{\partial}{\partial y} \big)}_{= \tilde{ \mathcal{L}}_2^{(y)}}  +
\underbrace{\frac{\partial}{\partial z} \big( a(x,y,z) \frac{\partial}{\partial z} \big)}_{= \tilde{ \mathcal{L}}_3^{(z)}} \Big) u(x,y,z).
\]
Applying the differential operator $\tilde{\mathcal{L}}_1^{(x)}$ to a polynomial $u \in \mathbb{P}_{n,n,n}$ can be written analogously to~\eqref{eq:CoeffTensorEquation3D} as
\[
\mathcal{V} = \mathcal{U} \times_1 \underbrace{S_1 D_1 S_{\rev{0}}^{-1} M^{(1)}[\mathbf{a}_1] D_1}_{= L_1^{(x)}} \times_2 \underbrace{M^{(2)}[\mathbf{a}_2] S_1 S_0}_{= L_1^{(y)}} \times_3 \underbrace{M^{(2)} [ \mathbf{a}_3] S_1 S_0}_{= L_1^{(y)}}
\]
where $\textbf{a}_1 \in \R^{n_1\rev{+1}}, \textbf{a}_2 \in \R^{n_2\rev{+1}},\textbf{a}_3 \in \R^{n_3\rev{+1}}$ denote the coefficients for univariate Chebyshev interpolation of $a_1,a_2,a_3$ as in Section~\ref{sec:1dNonConst}. The multiplication matrices $M[\textbf{a}]$, differentiation matrices $D_1$ and transformation matrices $S_0,S_1$ are defined as in Section~\ref{sec:Discretization}. We obtain analogous discretizations for $\tilde{\mathcal{L}}_2^{(y)}$ and $\tilde{\mathcal{L}}_3^{(z)}$. Observe that the discretization has the same symmetric structure as the CP decomposition~\eqref{eq:SpecialCPformForLaplacelike}. Thus, we can transform~\eqref{eq:3DReducedSystem} to a Laplace-like equation~\eqref{eq:LaplaceLikeEquation} by defining $U,V,W,\Breve{\mathcal{F}}$ as in~\eqref{eq:DefinitionMatricesLaplaceLike} and use the recursive solver. 

\subsubsection{Diffusion problems with higher rank coefficient} \label{sec:DiffusionPreconditioning}
Most diffusion problems~\eqref{eq:DiffusionProblem} arising in the study of groundwater flow and uncertainty quantification~\cite{Cohen10,Garreis17,Lord14,Schwenck15,Todor07,Ullmann12} do not have a rank-$1$ coefficient $a(x,y,z)$. The coefficient is often given by a truncated expansion as \rev{a} sum of separable functions of the form 
\[
a(x,y,z) = \sum_{r=1}^{R_a} a_{1}^{(r)}(x) a_2^{(r)}(y) a_3^{(r)}(z),
\]
with $R_a > 1$ and univariate functions $a_1^{(r)},a_2^{(r)},a_3^{(r)} :[-1,1] \to \R$.
The PDE~\eqref{eq:DiffusionProblem} can now be written as
\begin{equation}\label{eq:DiffusionHighRankPDE}
- \sum_{r=1}^{R_a} \bigtriangledown \cdot ( a_{1}^{(r)}(x) a_2^{(r)}(y) a_3^{(r)}(z) \bigtriangledown u(x,y,z)) = f(x,y,z).
\end{equation}
Following the ideas in Section~\ref{sec:SeparableDiffusion}, we can discretize~\eqref{eq:DiffusionProblem} for each separable function $a_{1}^{(r)}(x) a_2^{(r)}(y) a_3^{(r)}(z)$. Adding these discretizations yields a discretization of the form~\eqref{eq:DiscretizedPDE} with $R=3R_a$. Since $R > 3$, we can not find a Laplace-like formulation of~\eqref{eq:3DReducedSystem} and we can not use the recursive solver.

We can, however, use preconditioned GMRES~\cite{Golub13} to compute solutions of~\eqref{eq:Structure3DReducedSystem} seen as \rev{tensor-valued} linear system. 
Throughout this work\rev{,} we restart GMRES every 15 iterations. As preconditoner we employ the recursive solver to solve Equation~\eqref{eq:3DReducedSystem} for a discretization of the same diffusion problem~\eqref{eq:DiffusionHighRankPDE} with the coefficient $a(x,y,z)$ replaced by a separable coefficient $b(x,y,z)$. 
\rev{This can be seen as effectively solving the system 
\begin{align*}
\Big(\sum_{r=1}^3 \hat{L}_r^{(z)}[b] \otimes \hat{L}_r^{(y)}[b] \otimes \hat{L}_r^{(x)}[b]\Big)^{-1} \Big(\sum_{r=1}^R \hat{L}_r^{(z)}[a] \otimes \hat{L}_r^{(y)}[a] \otimes \hat{L}_r^{(x)}[a]\Big) \textrm{vec}(\rev{\mathcal{U}^{222}}) = \\
\Big(\sum_{r=1}^3 \hat{L}_r^{(z)}[b] \otimes \hat{L}_r^{(y)}[b] \otimes \hat{L}_r^{(x)}[b]\Big)^{-1} \textrm{vec}(\hat{\mathcal{F}}),
\end{align*}
where $\hat{L}_r^{(x)}[a],\hat{L}_r^{(y)}[a],\hat{L}_r^{(z)}[a]$ and $\hat{L}_r^{(x)}[b],\hat{L}_r^{(y)}[b],\hat{L}_r^{(z)}[b]$ denote the matrices in~\eqref{eq:Structure3DReducedSystem} based on discretizations of the PDE with coefficient $a$ and $b$ respectively. The application of the inverse of $\sum_{r=1}^3 \hat{L}_r^{(z)}[b] \otimes \hat{L}_r^{(y)}[b] \otimes \hat{L}_r^{(x)}[b]$ can be computed by solving a Laplace-like equation~\eqref{eq:LaplaceLikeEquation}. }

\begin{figure}
    \centering
    \includegraphics[width=0.44\textwidth]{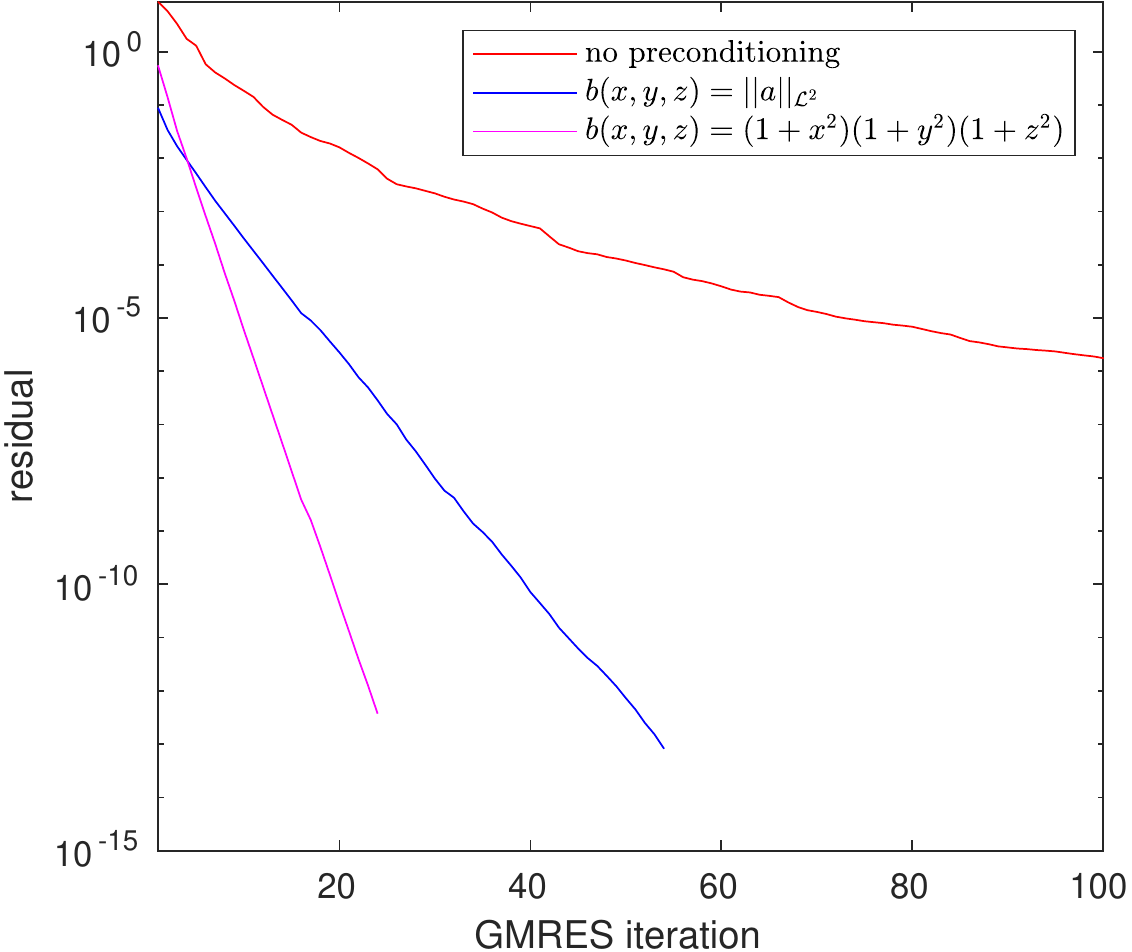}
    \caption{Convergence rate of preconditioned GMRES for solving Equation~\eqref{eq:3DReducedSystem} for the Diffusion problem~\eqref{eq:DiffusionHighRankPDE} with coefficient $a(x,y,z) = (1+x^2)(1+y^2)(1+z^2) + \exp(x+y+z)$. We compare no preconditioning and preconditioning with a separable coefficient $b(x,y,z)$ as described in Section~\ref{sec:DiffusionPreconditioning}.  }
    \label{fig:Preconditioning}
\end{figure}

From now on, we consider the rank-$2$ coefficient $a(x,y,z) = (1+x^2)(1+y^2)(1+z^2) + \exp(x+y+z)$. We use our global spectral method with $n = n_1 = n_2 = n_3 = 30$ to solve the Diffusion problem~\eqref{eq:DiffusionHighRankPDE} with known solution $u^*(x,y,z) = \sin(\pi x)  \sin(\pi y) \sin(\pi z)$, \rev{from which we explicitly compute the right hand side $f$ and Dirichlet boundary conditions}. In Figure~\ref{fig:Preconditioning}, we display the convergence rates for GMRES with preconditioners based on the constant coefficient $b_1(x,y,z) = ||a||_{\mathcal{L}^2}$ and the separable coefficient $b_2(x,y,z) = (1+x^2)(1+y^2)(1+z^2)$. Both preconditioners yield solutions $u$, for which the error $||u-u^*||_\infty$ is close to machine precision. Solving with $b_1$ takes \rev{$4.87$} seconds. For separable coefficient $b_2$ fewer iterations are necessary and the computation only takes \rev{$2.14$} seconds. For comparison, reshaping as in Section~\ref{sec:RuntimeComparison} would take \rev{$107$} seconds. 

\subsubsection{Helmholtz equation with non-constant coefficients} \label{sec:HelmholrtNonConst}
\rev{Next}, we study a variable coefficient Helmholtz equation as in~\cite[Example 2]{Townsend15}. We consider the PDE 
\begin{equation}\label{eq:HelmholtzNonConst}
        \Delta u(x,y,z) + \kappa(x,y,z) u(x,y,z) = f(x,y,z),
\end{equation}
with function $\kappa(x,y,z) = \sqrt{x+y+z+42}$. For $n = n_1=n_2=n_3 = 30$, we compute a discretization~\eqref{eq:DiscretizedPDE} by computing \rev{an approximate} CP decomposition of $\mathcal{A}$ \rev{with $R = 10$} as defined as in Section~\ref{sec:Discretize3DNonConst} using \rev{tensorlab~\cite{tensorlab}}. \rev{In $2.08$ seconds, we obtain a CP decomposition with error $2.04\cdot10^{-9}$ in the uniform norm.}
This leads to a discretized PDE, for which we solve Equation~\eqref{eq:3DReducedSystem} with preconditioned GMRES with restarting as in Section~\eqref{sec:DiffusionPreconditioning}. In the preconditioner we solve Equation~\eqref{eq:3DReducedSystem} for the Helmholtz equation~\eqref{eq:HelmholtzEquati} with constant coefficient $||\kappa||_{\mathcal{L}^2}$. 

When the right hand side and Dirichlet boundary conditions are chosen to match the known solution $u^*(x,y,z) = \sin(\pi x) \sin(\pi y) \sin(\pi z)$, solving with preconditioned GMRES takes $\rev{1.04}$ seconds. 
The computed solution $u$ satisfies $||u-u^*||_\infty \approx \rev{2.38\cdot 10^{-10}}$, where we estimate the uniform norm using $1\,000$ sample points.

\textit{Remark.} \rev{
To obtain very high accuracy in the PDE solution, we would need a very accurate CP decomposition. 
The efficient computation of accurate CP decompositions is still subject to research~\cite{Singh21}. 
Here, we could avoid computing a CP decompositoin of $\mathcal{A} \in \R^{\rev{93\times 93 \times 93}}$ by instead  approximating the operators $\Delta$ and $\kappa(x,y,z) \mathcal I$ separately, where where $\mathcal{I}$ denotes the identity operator. 
Discretizing $\kappa(x,y,z) \mathcal I$ only involves a CP decomposition of a tensor in $\R^{\rev{31 \times 31 \times 31}}$, which can be computed for $R=7$ in only $\rev{0.09}$ seconds with error $1.93\cdot 10^{-11}$ in the uniform norm. 
The two resulting discretizations in the form of~\eqref{eq:3dRawDiscretization} can be added to obtain a discretization of $\Delta + \kappa(x,y,z) \mathcal I$.}

\subsubsection{\rev{Helmholtz equation with unknown solution}}\label{sec:UnknownSolution}

\begin{figure}
    \centering
\subfloat{{\includegraphics[width=0.44\textwidth]{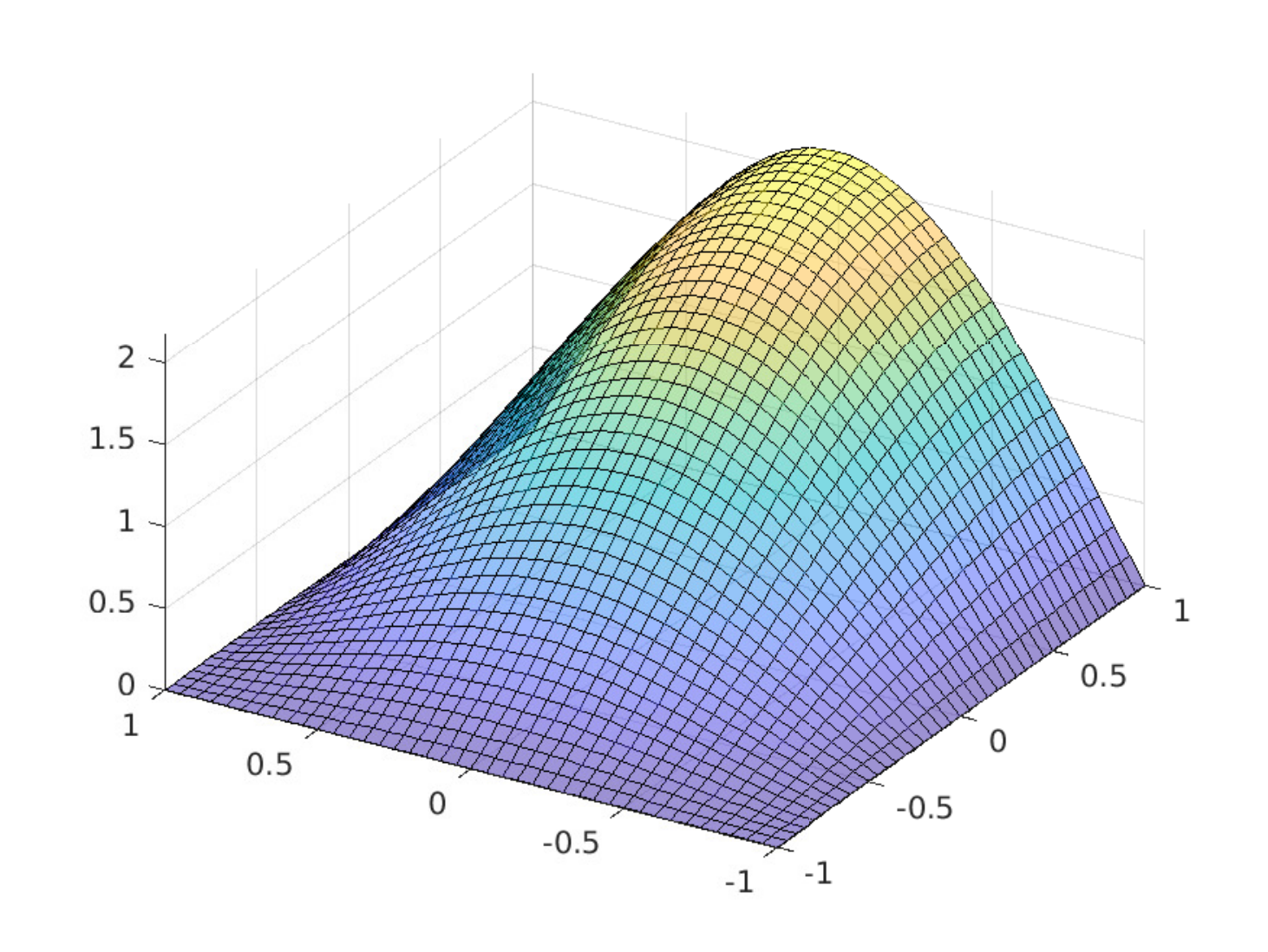}}}%
\qquad
\subfloat{{\includegraphics[width=0.44\textwidth]{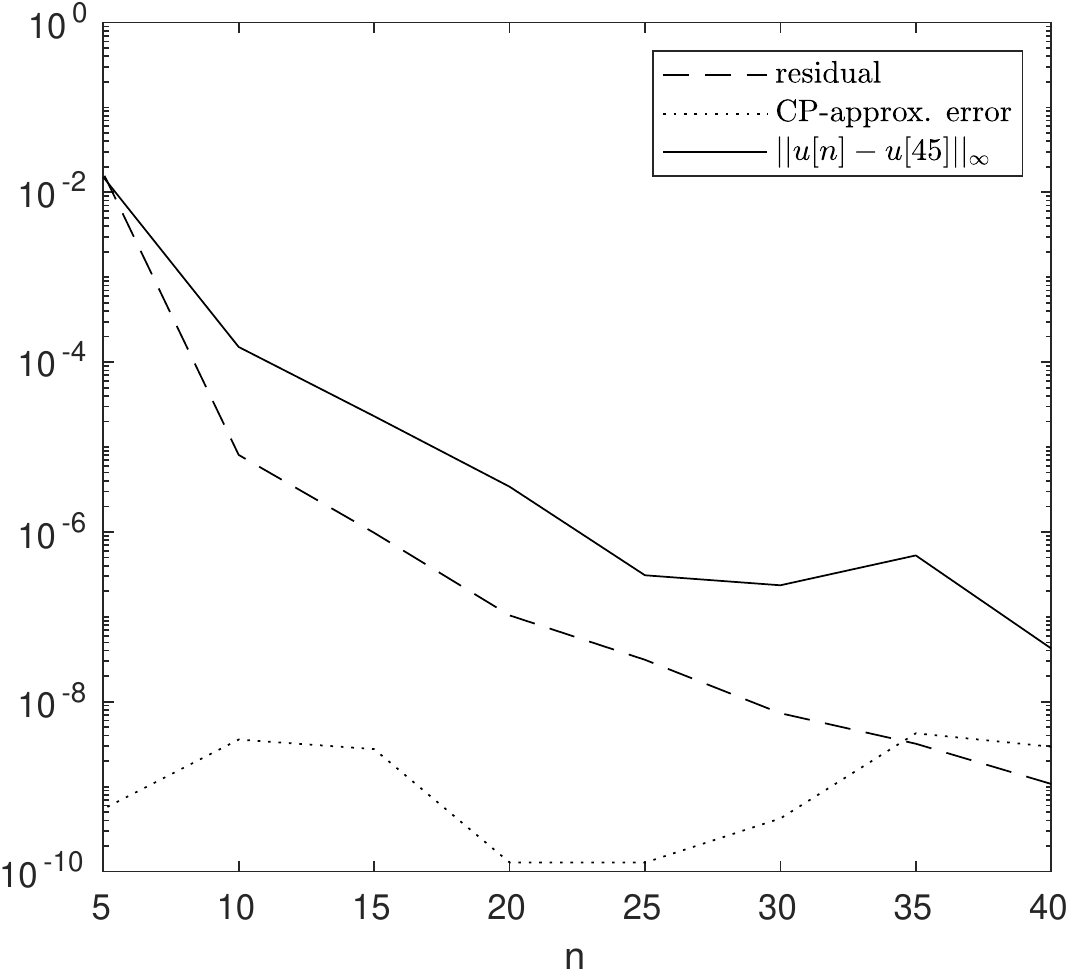}}}%
    \caption{\rev{ 
    Left: Let $u[n] \in \mathbb{P}_{n,n,n}$ denote the solution computed for the Helmholtz equation described in Section~\ref{sec:UnknownSolution}.   We plot $u[45](x,y,1/4)$. Right: For various $n$ we plot $||u[n]-u[45]||_{\infty}$ estimated at $1\,000$ sample points. Additionally, we plot the combined residual defined as maximum of the residual of~\eqref{eq:DiscretizedPDE} and the residuals of~\eqref{eq:GeneralDiscretizedBoundaryConditions} in the uniform norm. Further, we plot the error of the CP approximation of $\mathcal{A}$ in the uniform norm.}}
    \label{fig:UnknownSolution}
\end{figure}

\rev{As a final stationary problem, we consider the Helmholtz equation~\eqref{eq:HelmholtzNonConst} with non-constant coefficient $\kappa(x,y,z) = \sqrt{x+y+z+42}$ and $f(x,y,z) = 1$.
We use mixed boundary conditions with zero Neumann boundary conditions on the right and zero Dirichlet boundary conditions on all other sides as in Section~\ref{sec:Neumann}.
As in Section~\ref{sec:HelmholrtNonConst}, we compute a CP decomposition of $\mathcal{A}$ and solve using preconditioned GMRES. 
In Figure~\ref{fig:UnknownSolution}, we display the computed solution and the error decay for different polynomial degrees. 
We observe that the residual decays when the polynomial degree is increased, which indicates that the solutions become more accurate.
}

\subsection{Time-dependent problems}\label{sec:TimeDependent}
In this section, we introduce an implicit Euler scheme to solve parabolic PDEs of the form \[\frac{\partial}{\partial t} u(x,y,z,t) + \L u(x,y,z,t) = 0,\]
where $\L$ is a linear partial differential operator acting only on the spatial variables $x,y,z$. The system is complemented with boundary conditions. We are interested in the time evolution starting from a given initial function $u(x,y,z,0) = u_0(x,y,z)$. We discretize the equation in time using the uniform step length $h$. For each $\tau = 1,2,\dots$, we compute $u_\tau(x,y,z) \approx u(x,y,z,\tau h)$ as solution of the stationary linear PDE
\begin{equation}\label{eq:ImplicitEulerAbstract}
    (\mathcal{I} - h \L) u_{\tau+1} = u_\tau,
\end{equation}
We approximate each function $u_\tau$ by a polynomial of the form~\ref{eq:ChebyshevFormat} represented by the coefficient tensor $\mathcal{U}_\tau$. The initial $\mathcal{U}_0$ is computed via \rev{tensorized} Chebyshev interpolation. For $\tau=1,2,\dots$, we obtain $\mathcal{U}_\tau$ by applying our spectral method to solve Equation~\eqref{eq:ImplicitEulerAbstract}. We discretize the operator $(\mathcal{I} - h \L)$ directly like the Helmholtz equation in Section~\ref{sec:Helmholtz}. This allows us to employ the recursive solver. Note that the right hand side is represented in terms of an \rev{ultraspherical} basis. Hence, the computation of $\mathcal{U}_\tau$ requires multiplying $\mathcal{U}_{\tau-1}$ with appropriate basis transformation matrices.

We demonstrate this implicit Euler scheme for the parabolic PDE studied in~\cite[Section 6.1]{Petersdorff04}. The function $u^*(x,y,z,t) = \exp\rev (-3 \pi^2 t\rev ) \sin(\pi x) \sin(\pi y) \sin(\pi z)$ satisfies the parabolic PDE \[
\frac{\partial}{\partial t} u + \Delta u = 0,
\] on the domain $[-1,1]^3$ with homogeneous Dirichlet boundary conditions. We use an implicit Euler scheme for $u_0(x,y,z) =  \sin(\pi x) \sin(\pi y) \sin(\pi z)$ and compare  our global spectral method  to a finite difference method. 
In each timestep of the finite difference method, we solve a linear system with a sparse Kronecker-structured finite difference matrix using the backslash operator in MATLAB. 

\begin{figure}[ht]
    \centering
\subfloat[Results for $n=20$ and $h = 10^{-2}$.]{{\includegraphics[width=0.44\textwidth]{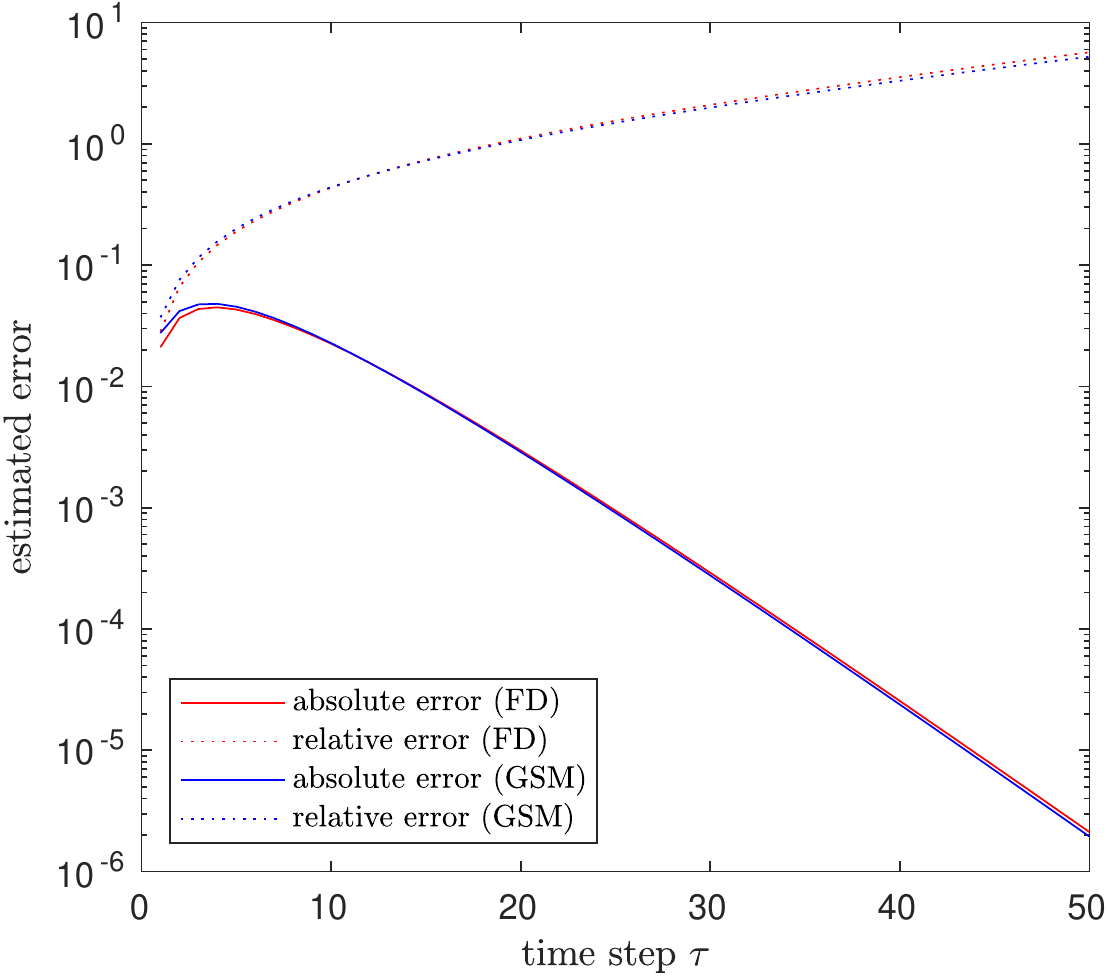}}}%
\qquad
\subfloat[Results for $n=20$ and $h = 10^{-4}$.]{{\includegraphics[width=0.44\textwidth]{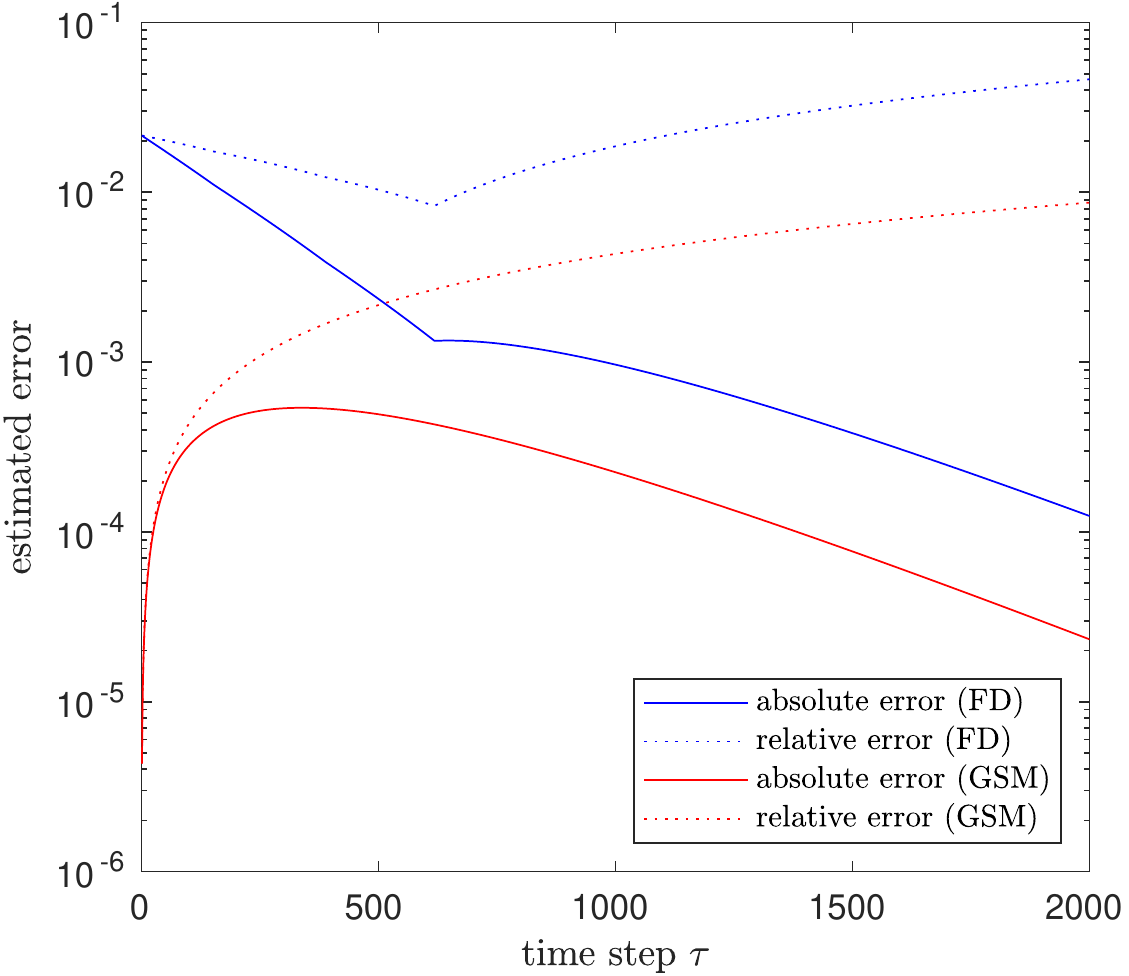}}}%
    \caption{
    Error plots for the parabolic PDE example in Section~\ref{sec:TimeDependent}. We apply an implicit Euler scheme to compute approximations $u_\tau(x,y,z) \approx u(x,y,z,\tau h)$ and compare the solutions for  the global spectral method (GSM) with solutions in $\mathbb{P}_{n,n,n}$ and a finite difference discretization (FD) on a regular $(n \rev{+1)} \times (n \rev{+1)} \times (n \rev{+1)}$ grid. The absolute error $||u^*(\rev{\cdot},\rev{\cdot},\rev{\cdot},\tau h)-u_\tau(\rev{\cdot},\rev{\cdot},\rev{\cdot})||_\infty$ is estimated in each timestep from evaluations at $100$ sample points. The relative error is computed from the estimation of the absolute error.}
    \label{fig:TimeDependent}
\end{figure}

The time evolution of the errors is displayed in Figure~\ref{fig:TimeDependent}. We observe that for $n=20$ and $h=10^{-2}$ both approximations lead to very similar errors, but the computation time for $50$ implicit Euler steps decreases from  $\rev{6.22}$ seconds for the finite difference method to $\rev{0.93}$ seconds for our global spectral method. In this case the error is dominated by the implicit Euler scheme for both methods. In contrast, for $n=30$ and $h = 10^{-4}$, the errors for the spectral method are smaller. For large $\tau$ the error of both approaches is dominated by the time discretization via the implicit Euler scheme. However, in the initial $750$ time steps the error of the spatial discretization dominates for the finite difference scheme, whereas the global spectral method is able to represent $u_\tau$ accurately. 

\textit{Remark.} The methods presented in this work can also be used to solve two-dimensional parabolic PDEs on rectangles by treating time as third space variable in the discretization of the operator. We refer to~\cite{Townsend15} for more details on this approach.

\subsection{Eigenvalue problems}

The methods presented in this work can be extended to solve PDE eigenvalue problems, in which we search eigenvalues $\lambda$ and eigenfunctions $u$ satisfying the equation $\L u = \lambda u$ complemented with homogeneous Dirichlet boundary conditions. We are particularly interested in finding the eigenvalue with minimal absolute value. For this purpose we employ the inverse iteration algorithm~\cite{Golub13}. Starting from an initial function $u_0$ we iteratively compute an approximation of the eigenfunction. For ${\rev{s}}=1,2,\dots$, we compute $u_{{\rev{s}}}$ as solution of the PDE \begin{equation}\label{eq:InverseIterationStep}
    \L u_{\rev{s}} =  \frac{u_{{\rev{s}}-1}}{||u_{{\rev{s}}-1}||_{\mathcal L^2}}.
\end{equation} We approximate the eigenvalue using the Rayleigh quotient $\frac{1}{\lambda} \approx \frac{\langle u_{{\rev{s}}-1}, u_{\rev{s}} \rangle}{\langle u_{{\rev{s}}-1} , u_{{\rev{s}}-1} \rangle}$, where $\langle\cdot , \cdot \rangle$ denotes the standard $\mathcal{L}^2$ scalar product.
We again proceed by discretizing the differential operator and the function $u$ to solve Equation~\eqref{eq:InverseIterationStep} using the spectral method introduced in Section~\ref{sec:SpectralMethod}.

Let the functions $u,v \in \mathbb{P}_{n_1,n_2,n_3}$ be given in the form of~\eqref{eq:ChebyshevFormat}. We evaluate the norm as $||u||_{\mathcal{L}^2}^{\rev{2}} = \langle u , u \rangle$  and the scalar product $\langle u , v \rangle$ by interpolating the function $uv$ using \rev{tensorized} Chebyshev polynomial basis functions as in~\cite{Stroessner21,Mason02}. In a second step, we integrate the approximation of $uv$ using the exact values for the integrals of the basis functions~\cite[Theorem 19.2]{Trefethen13}.

\begin{figure}[ht]
    \centering
    \includegraphics[width=0.44\textwidth]{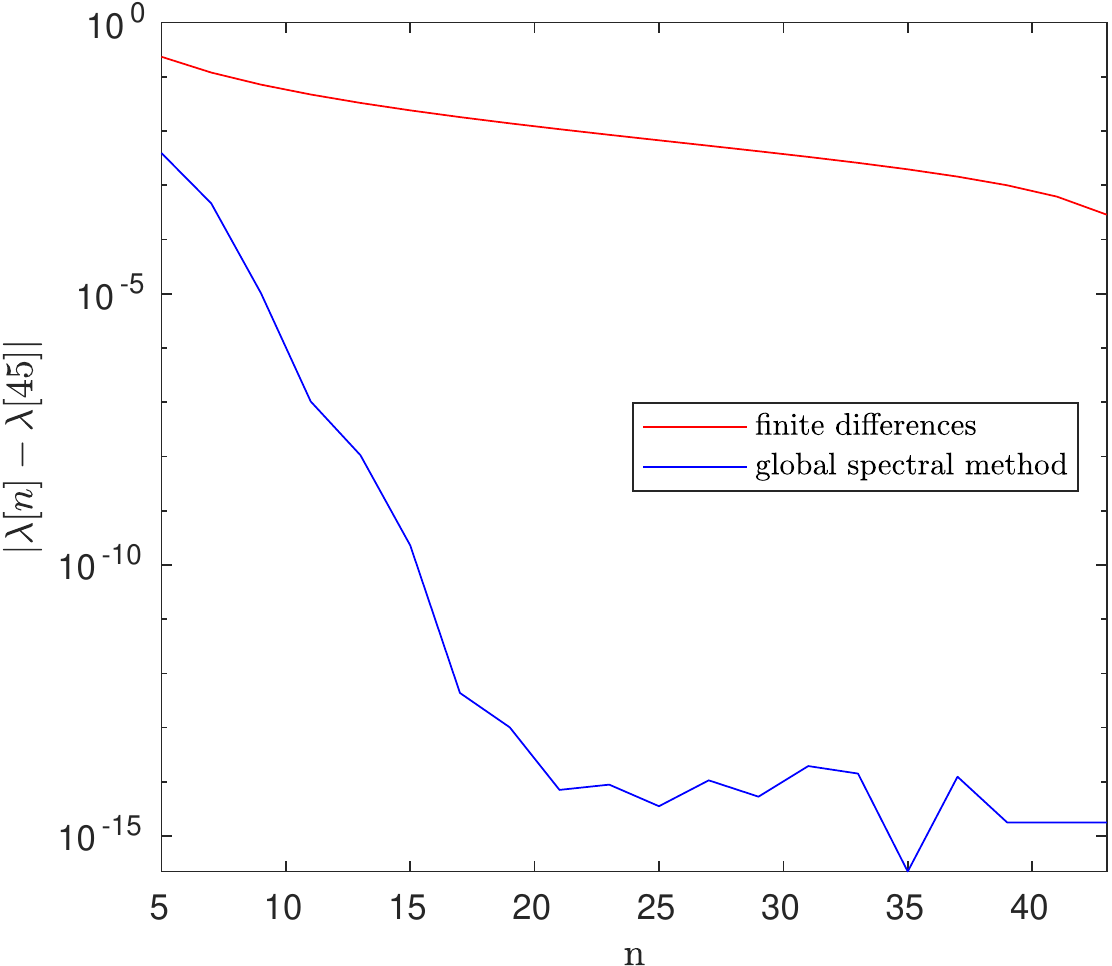}
    \caption{Eigenvalue convergence rate for the PDE eigenvalue problem~\eqref{eq:ellipticPDEeigenproblem} computed as Rayleigh coefficients from $u_{50}$ in the inverse iteration method, where the solutions $u_{\rev{s}}$ are computed via the global spectral method with $u_{\rev{s}} \in \mathbb{P}_{n,n,n}$ and a finite difference method on a regular $(n\rev{+1)} \times (n\rev{+1)} \times (n\rev{+1)}$ grid. We denote the eigenvalue obtained for $n$ by $\lambda[n]$. \rev{Since the exact $\lambda$ is not known, we plot $|\lambda[n]-\lambda[45]|$ for both algorithms, where $\lambda[45]$ is computed with the respective algorithm.}}
    \label{fig:Eigenvalues}
\end{figure}

We test our method for the elliptic PDE eigenvalue problem
\begin{equation}\label{eq:ellipticPDEeigenproblem}
    -\Delta u(x,y,z) + v(x,y,z) u(x,y,z) = \lambda u(x,y,z),
\end{equation}
with potential $v(x,y,z) = \sin(\pi/2(x+1))  \sin(\pi/2 (y+1))  \sin(\pi/2(z+1))$ as in~\cite{Kressner12,Hackbusch12}. 
Due to the potential $v(x,y,z)$ the inverse iteration steps~\eqref{eq:InverseIterationStep} are Helmholtz equations~\eqref{eq:HelmholtzEquati} with non-constant, separable coefficients. 
We discretize these by adding a discretization  of $v(x,y,z) \mathcal{I}$ \rev{with $R=1$} to a discretization of the Laplacian. 
The resulting discretized PDE~\eqref{eq:DiscretizedPDE} has $R=4$ and we apply preconditioned GMRES with restarting using a Helmholtz equation with constant coefficient $||v||_{\mathcal{L}^2}$ to compute solutions as in Section~\ref{sec:HelmholrtNonConst}. 
We compare our global spectral method to a finite difference scheme on a regular grid for solving~\eqref{eq:InverseIterationStep}.
In Figure~\ref{fig:Eigenvalues} we observe that the eigenvalues computed as Rayleigh coefficients converge at a faster rate for the global spectral method \rev{than} the finite difference approach. Already for $n=20$ the global spectral method reaches an error close to machine precision. 

\textit{Remark.} We would like to point out that our approach can be used to compute a basis for the (continuous) Krylov subspace, which is used in the Arnoldi method~\cite{Arnoldi51,Golub13} and in the Least-squares spectral method in~\cite{Hashemi21}.

\section{Conclusions}
In this work, we derive a global spectral method for solving three-dimensional linear PDEs on cubes with very high accuracy. 
We demonstrate that the Laplace-like equations arising for certain PDEs can be solved efficiently with the blocked recursive solver~\cite{Chen19}. 
Our numerical experiments show that applying this solver directly or as \rev{a} preconditioner vastly outperforms all existing methods. 
The versatility of our method is presented by the extension to eigenvalue and time-dependent problems.

\paragraph{Future work.}
The computational complexity of our global spectral method is heavily  influenced by the storage needed to store the full coefficient tensors $\mathcal{U}$ and $\mathcal{F}$ for representing the solution and right-hand side, respectively. 
For certain problems these tensors admit good low-rank approximations. Applying such approximations leads to so\rev{-}called functional low-rank approximations~\cite{Bigoni16,Stroessner21,Gorodetsky19,Hashemi17,Soley21}, which have the potential to drastically reduce the storage complexity. 
Exploiting this potential requires determining suitable approximation formats and using a specialized solver for the chosen formats.  
For instance, in~\cite{Shi21} it is shown that under certain conditions the solution of Laplace-like equations~\eqref{eq:LaplaceLikeEquation} can be represented in tensor train or Tucker format when the right hand side is given in the same format.
The functional tensor train format~\cite{Bigoni16,Chertkov21b,Gorodetsky19,Soley21}, could potentially be used to extend the global spectral method to higher dimensional, linear PDEs on hypercubes.
The computational complexity of time-dependent problems could be reduced using rank-adaptive, dynamical low-rank approximations~\cite{Bachmayr20,Dektor20,Dektor19,Koch07}.

\rev{The current implementation of our method is far from Chebfun-like~\cite{Driscoll14}. Many more ingredients, such as adaptivity and input parsing, are necessary to reach the level of Chebfun. These modifications are subject to future work.}

\end{document}